\documentclass[11pt,a4paper,twoside,reqno]{amsart}

\usepackage[english]{babel}
\usepackage{amsmath}
\usepackage{amssymb} 
\usepackage{amsfonts}
\usepackage[all,color]{xy}
\usepackage{hyperref}
\usepackage{lscape}
\usepackage{graphicx}
\usepackage{float}

 \numberwithin{equation}{section}
 \usepackage{color}
  \usepackage{hyperref}
\usepackage{enumerate,xspace}
\usepackage{graphicx}

\renewcommand{\tilde}{\widetilde}

\renewcommand{\phi}{\varphi}
\renewcommand{\epsilon}{\varepsilon}

 \newtheorem{proposition}{Proposition}[section] 
 \newtheorem{lemma}[proposition]{Lemma}
 \newtheorem*{lemma*}{Lemma}
 \newtheorem{theorem}[proposition]{Theorem}
 \newtheorem{corollary}[proposition]{Corollary}
 \theoremstyle{definition}
 \newtheorem{definition}[proposition]{Definition}
 \newtheorem*{definition*}{Definition}
 \newtheorem{example}[proposition]{Example}
 \newtheorem*{example*}{Example}
 \newtheorem{remark}[proposition]{Remark}
 
\numberwithin{equation}{section}

\newcommand{\ca}{\ensuremath{\mathcal A}\xspace}

\newcommand{\cc}{\ensuremath{\mathcal C}\xspace}
\newcommand{\cd}{\ensuremath{\mathcal D}\xspace}

\newcommand{\cl}{\ensuremath{\mathcal L}\xspace}

\newcommand{\cn}{\ensuremath{\mathcal N}\xspace}
\newcommand{\cp}{\ensuremath{\mathcal P}\xspace}
\newcommand{\cq}{\ensuremath{\mathcal Q}\xspace}
\newcommand{\crr}{\ensuremath{\mathcal R}\xspace}
\newcommand{\cs}{\ensuremath{\mathcal S}\xspace}
\newcommand{\ct}{\ensuremath{\mathcal T}\xspace}

\newcommand{\bbn}{\ensuremath{\mathbb N}\xspace}

\newcommand{\bs}{\ensuremath{\mathbb S}\xspace}

\newcommand{\Fsk}{\textnormal{\bf Fsk}\xspace}

\newcommand{\Tam}{\mathrm{Tam}}

\DeclareMathOperator{\End}{End}
\DeclareMathOperator{\im}{im}

\newcommand{\LNC}{\textnormal{LBC}}

\newcommand{\myt}{t}
\newcommand{\myl}{\ell}
\newcommand{\mylam}{\lambda}

\def\ox{\otimes}
\def\x{\times}

\renewcommand{\>}{\rangle}

\newcommand{\ps}{\ensuremath{{}_{\textrm{ps}}}\xspace}
\newcommand{\wk}{\ensuremath{{}_{\textrm{w}}}\xspace}

\newcommand{\op}{\ensuremath{{}^{\textrm{op}}}\xspace}

\newcommand{\Cat}{\ensuremath{\mathbf{Cat}}\xspace}
\newcommand{\Ord}{\ensuremath{\mathbf{Ord}}\xspace}
\newcommand{\Mod}{\ensuremath{\mathbf{Mod}}\xspace}

\newcommand{\Algs}{\ensuremath{\textnormal{-Alg}_{\textnormal{s}}}\xspace}
\newcommand{\Algl}{\ensuremath{\textnormal{-Alg}_\ell}\xspace}
\newcommand{\Algw}{\ensuremath{\textnormal{-Alg}_{\textnormal{w}}}\xspace}

\newcommand{\Skews}{\ensuremath{\mathbf{Skew}_{\textnormal{s}}}\xspace}
\newcommand{\Skeww}{\ensuremath{\mathbf{Skew}_{\textnormal{w}}}\xspace}

\newcommand{\nc}{\textnormal{nColax-}}
\newcommand{\clx}{\textnormal{Colax-}}

\newcommand{\ord}[1]{\ensuremath{\mathbf{#1}}}

\definecolor{lightgrey}{rgb}{0.666666,0.666666,0.666666}
\definecolor{pinegreen}{rgb}{0.15,0.7,0.15}

\begin{document}

\title{Free skew monoidal categories}
\author{John Bourke, Stephen Lack}
\address{Department of Mathematics, Macquarie University NSW 2109, 
Australia}
\email{john.d.bourke@mq.edu.au}
\address{Department of Mathematics, Macquarie University NSW 2109, 
Australia}
\email{steve.lack@mq.edu.au}
\begin{abstract}
In the paper {\em Triangulations, orientals, and skew monoidal categories}, the free
monoidal category \Fsk on a single generating object was described. We
sharpen this by giving a completely explicit description of \Fsk, and so of the free skew monoidal category on any
category. 
As an application we describe adjunctions between the operad for skew
monoidal categories and various simpler operads.  
For a particular such operad $\cl$, we identify skew monoidal categories with certain colax $\cl$-algebras.
\end{abstract} 
\date\today
\maketitle

\section{Introduction}

A skew monoidal category is a category \cc equipped with a functor
$\cc^2\to \cc$ whose effect on objects we write as $(a,b)\mapsto ab$,
an object $i\in\cc$, and natural transformations
\[ \xymatrix @R0pc { 
(ab)c \ar[r]^{\alpha} & a(bc) \\
ia \ar[r]^{\lambda} & a \\
a \ar[r]^{\rho} & ai
} \]
satisfying five coherence conditions. When the maps $\alpha$, $\rho$,
and $\lambda$ are invertible, we recover the usual notion of monoidal
category. 

While this might seem like a mindless generalisation, it turns out
that there are important examples of skew monoidal categories which
are not monoidal. The first such class of examples arises from
quantum algebra, and is due to
Szlach\'anyi \cite{Szlachanyi-skew}: he realised that bialgebroids can be described using
skew monoidal categories. Specifically, a bialgebroid with base ring
$R$ is the same thing as a skew monoidal closed structure on the
category $R$-$\Mod$  of $R$-modules. 

A second class of examples arises from the intersection of homotopical algebra
and 2-category theory: a host of naturally occurring skew monoidal closed structures
on Quillen model categories that arise in 2-dimensional universal
algebra were described in \cite{bourko-skew}. 
These examples are monoidal in a homotopical sense, in that they yield
genuine monoidal closed structures on the associated homotopy categories.

Unlike the situation for monoidal categories, it is not the case for
skew monoidal categories that all diagrams built up out of the
structure maps commute: for example, the composite 
\[ \xymatrix{ ii \ar[r]^{\lambda} & i \ar[r]^{\rho} & ii } \]
is not the identity, and so the ``coherence problem'' for skew
monoidal categories is not a trivial one. One way to formulate this
coherence problem is to ask what is the {\em free skew monoidal
  category} on a given category. An answer to this question was given
in \cite{skewcoherence}.

As was observed in \cite{skewcoherence}, the structure of a skew
  monoidal category  is {\em clubbable}, in the sense of
  \cite{Kelly:clubs}; equivalently, it can be given in terms of a
  {\em plain operad in \Cat}, where by ``plain'', we mean that there are
  no actions of the symmetric groups. It then follows that in order to
  describe the free skew monoidal category on a general category \cc
  it suffices to do it on the terminal category, and in fact this is
  what is done in \cite{skewcoherence}.

The free skew monoidal category on an object, called \Fsk in
\cite{skewcoherence}, is determined by the following universal
property. There is a designated object $X\in\Fsk$ (``the generator'') and for any
skew monoidal category \cc, evaluation at $X$ determines a bijection
between the set of (strict) monoidal functors from $\Fsk$ to $\cc$ and
objects of $\cc$.

An example of a skew monoidal category is the category $\Ord_\bot$ of
finite non-empty ordinals, with morphisms the functions which preserve
both order and bottom element. The product is given by ordinal sum,
and the unit object is the ordinal $\ord 1=\{0\}$. This is strictly associative,
but the maps $\lambda$ and $\rho$ are non-invertible. By its universal
property, the free skew monoidal category \Fsk on one object has a unique
structure-preserving functor to $\Ord_\bot$ which sends the generator
to $\ord 1$. A key result of \cite{skewcoherence} was that this functor is
faithful, so that the morphisms of \Fsk can be represented as certain
functions between finite sets. 

While the objects of \Fsk were described in an entirely explicit way,
the morphisms were not. The main goal of
this paper is to remedy this, by giving a completely explicit
condition characterising the morphisms. 

As an application, we construct various adjunction between the operad
\cs for skew monoidal categories and various simpler operads $\ct$.  In each
case we have an operad map $F\colon\cs \to \ct$ and we show that $F$ has a 
left or right adjoint in each component.  By the usual ``doctrinal adjunction"
results \cite{Kelly-doctrinal} this enables us to view skew monoidal categories
as colax/lax \ct-algebras. When \ct is the terminal operad, the unique
map $F\colon \cs\to\ct$ has both adjoints, and so any skew monoidal
category yields both a colax monoidal category and a lax monoidal
category. These processes lose structure, but choosing for \ct an only 
slighter more complex operad \cl, we find that colax \cl-algebras encode
the skew monoidal structure entirely.   These results are used in our
companion paper \cite{bourke-lack} which introduces and studies \emph{skew multicategories},
the multicategorical analogue of skew monoidal categories.

\subsection*{Acknowledgements}

Both authors acknowledge with gratitude the support of an Australian Research Council
Discovery Grant DP160101519.

\section{Background on clubs and operads}\label{sect:operads}

In this section we group together various facts about plain
\Cat-operads and clubs over \bbn. There is nothing particularly original
here, but we could not find any convenient reference containing
everything we need, which largely amounts to combining aspects of
\cite{Kelly:clubs}, \cite{Kelly-operads},  and \cite{Kelly-LaxAlg}.

Let \bbn denote the discrete category with objects the natural
numbers. The functor 2-category $[\bbn,\Cat]$ 
has a monoidal structure with tensor product $\circ$ given
by $(A\circ B)_n = \sum_{n=n_1+\ldots+n_k} A_k\x B_{n_1}\x\ldots\x
B_{n_k}$ and with unit $J$ given by $J_1=1$ and $J_n=0$ if $n\neq 1$.

A monoid in $[\bbn,\Cat]$ is called a {\em plain \Cat-operad}
\cite{Kelly-operads}; here the epithet ``plain'' serves to distinguish
these operads from the variant involving actions of the symmetric
groups. Since plain \Cat-operads are the only operads which appear in
this paper, we may sometimes simply call them operads.

We generally write $\mu\colon \ct\circ\ct\to\ct$ for the
multiplication and $\eta\colon J\to\ct$ for the unit of an operad
\ct. 
Explicitly, the components of the multiplication are ``substitution'' maps 
\[ \xymatrix @R0pc { \ct_k \x \ct_{n_1}\x\ldots\x \ct_{n_k} \ar[r] &
    \ct_{n_1+\ldots+n_k} \\
(g,f_1,\ldots,f_k) \ar@{|->}[r] & g(f_1,\ldots,f_k) }
\]
while the unit amounts to an object of $\ct_1$.

In the special case where all the $f$s are identities except
(possibly) for $f_i$, we sometimes write $g\circ_i f_i$ for
$g(f_1,\ldots,f_k)$. It is possible to reformulate the definition of
operad using only the operations $\circ_i$; for instance
$g(f_1,\ldots,f_k)$ can be constructed as 
\[ g(f_1,\ldots,f_k) =  (((g\circ_k f_k)\ldots\circ_2 f_2) \circ_1 f_1). \]

The 2-category $[\bbn,\Cat]$ is equivalent to the 2-category
$\Cat/\bbn$ of categories over \bbn, with the equivalence sending a
functor  $\ct\colon\bbn\to\Cat$ to the coproduct $\sum_{n\in\bbn} \ct_n\to\bbn$.
The monoidal structure on
$[\bbn,\Cat]$ can be transported across the equivalence to obtain a
monoidal structure on $\Cat/\bbn$ (although it was first defined
independently of that on $[\bbn,\Cat]$). A monoid in $\Cat/\bbn$ is
called a {\em club over $\bbn$} \cite{Kelly:clubs}; once again, in
this paper no other clubs are considered so we may simply call it a
club. 

There is a functor $E\colon \bbn\to\Cat$ sending the natural number
$n$ to the discrete category with $n$ objects. Left Kan extension
along $E$ determines a functor $[\bbn,\Cat]\to[\Cat,\Cat]$ which is
strong monoidal and so sends monoids to monoids; that is, it sends
plain \Cat-operads to 2-monads on \Cat. 

An algebra for an operad is an algebra for the corresponding 2-monad,
but these can also be described directly: a \ct-algebra is a category
$A$ equipped with functors $\ct_n\x A^n\to A$ for each $n$, satisfying
associativity and unit conditions. 

There is also another approach. For a category $A$ there is
an operad $\End(A)$ with $\End(A)_n=\Cat(A^n,A)$ and the substitution
maps for $\End(A)$ are given by actual substitution
\[ \Cat(A^k,A)\x \Cat(A^{n_1},A)\x\ldots\x \Cat(A^{n_k},A) \to
  \Cat(A^{n_1+\ldots+n_k},A). \]
For an operad $\ct$, an algebra structure on $A$ is the same as an
operad map $a\colon \ct\to\End(A)$.

Because of the 2-category
structure of $[\bbn,\Cat]$ it would be possible to consider various
weakening of the notion of operad. We do not do this, but we do
consider weak morphism of operads. Specifically, we consider colax
morphisms of operads; formally, these are analogous to opmonoidal
functors between monoidal categories. 

For operads \ct and \cq, then, a {\em colax morphism of operads} from
\ct to \cq is a morphism $F\colon \ct\to\cq$ in $[\bbn,\Cat]$ equipped
with 2-cells 
\[ \xymatrix{
\ct\circ\ct \ar[r]^{F\circ F}_{~}="2" \ar[d]_\mu & \cq\circ\cq
\ar[d]^{\mu}  & J \ar[d]_{\eta}  \ar@/^1pc/[dr]^{\eta}_(0.4){~}="4" \\
\ct \ar[r]_{F}^{~}="1" & \cq & 
\ct\ar[r]_{F}^{~}="3" & \cq 
\ar@{=>}"3";"4"_{\tilde{F}_0} 
\ar@{=>}"1";"2"_{\tilde{F}} 
} \]
satisfying one coassociativity and two counit condition, analogous to
those for opmonoidal functors. Such a colax morphism is said to be
{\em normal} if $\tilde{F}_0$ is an identity, and it is this normal
case in which we are primarily interested. Of course if
$\tilde{F}$ is also an identity we recover the usual (strict) notion of
morphism of operads. 

The components of $\tilde F$ have the form 
\[ \xymatrix{ F(g(f_1,\ldots,f_n)) \ar[r]^-{\tilde{F}} & F(g)
    ( F(f_1),\ldots,F(f_n)) } \] 
which need to be natural in $g$ and the $f_i$, as well as satisfying
the coassociativity and counit conditions. By coassociativity, these 
components can all be recovered from the special cases where all but one
of the $f_i$ is an identity, which then look like
\[ \xymatrix{ F(g\circ_i f) \ar[r]^-{\tilde{F}} & F(g)\circ_i
    F(f). } \]

If $F\colon \ct\to\cq$ and $G\colon \cq\to\cp$ are colax
morphisms of operads, then we may paste $\tilde{F}$ and
$\tilde{G}$ to give $G\circ F$ the structure of colax
morphism from \ct to \cp, normal if $F$ and $G$ are so. 

\begin{definition}
A {\em colax algebra} for an operad \ct is a category $A$
equipped with a colax morphism $\ct\to\End(A)$. The colax algebra is
{\em normal} if the colax morphism if so. 
\end{definition}

More explicitly, these involve functors $m\colon \ct_n\x A^n\to A$
for each $n$, and so in particular functors $m_x\colon A^n\to A$ for
each $n$ and each object $x\in \ct_n$. Then the colax structure
involves natural transformations 
\[ \xymatrix{
A^{k+n-1} \ar[d]_{A^{i-1}\x m_y\x A^{k-i}} 
\ar@/^1pc/[dr]^{m_{x\circ_i y}}_(0.42){~}="2" \\
A^k \ar[r]_{m_x}^{~}="1" & A 
\ar@{=>}"2";"1"_{\Gamma_{y,i,x}}  } \]
or in other words $\Gamma_{x,i,y}\colon m_{x\circ_i y} \to m_x\circ_i
m_y$.

Whether or not our algebras are strict or colax, we need to consider
various flavours of weak morphism between them. 

We defined $\End(A)$ so that $\End(A)_n=\Cat(A^n,A)$. More generally,
if $A$ and $B$ are categories, we may define $\<A,B\>\in[\bbn,\Cat]$
by $\<A,B\>_n=\Cat(A^n,B)$; thus $\<A,A\>$ is the underlying object in
$[\bbn,\Cat]$ of the operad $\End(A)$. This construction is functorial
in both arguments, covariant in $B$ and contravariant in $A$.

If now $f\colon A\to B$ is a functor, we may form the comma object 
\[ \xymatrix{
\{f,f\}_\ell \ar[r]^c_{~}="1" \ar[d]_d^{~} & 
\<B,B\> \ar[d]^{\<f,B\>}_{~} \\
\<A,A\> \ar[r]_{\<A,f\>}^{~}="2" & \<A,B\> 
\ar@{=>}"1";"2" } \]
in $[\bbn,\Cat]$ and $\{f,f\}_\ell$ inherits a unique operad structure for which $d$
and $c$ are morphisms of operads. An object of $(\{f,f\}_\ell)_n$
consists of functors $x\colon A^n\to A$ and $y\colon B^n\to B$, along
with a natural transformation $y.f^n\to f.x$.

\begin{definition}\label{defn:lax-morphism}
If $A$ and $B$ are colax \ct-algebras, with $F\colon
\ct\to\End(A)$ and $G\colon \ct\to\End(B)$ the corresponding 
colax morphisms of operads, a {\em lax morphism} from $A$ to $B$ is a
functor $f\colon A\to B$ equipped with a colax morphism
$H\colon \ct\to\{f,f\}_\ell$ for which $d H=F$ and $c H=G$.
\end{definition}

To give such an $H$ is to give natural transformations
\[ \xymatrix{
\ct_n\x A^n \ar[r]^{1\x f^n}_{~}="2" \ar[d] & \ct_n\x B^n \ar[d] \\
A \ar[r]_{f}^{~}="1" & B 
\ar@{=>}"2";"1"^{\tilde{f}}  } \]
satisfying associativity and unit conditions. 

In the case where $A$ and $B$ are strict algebras, so that $F$ and
$G$ are strict morphisms of operads, such an $H$ will itself be
a strict morphism of operads. This is not of course to say that $f$ is
strict as a morphism of algebras; its lax nature has been incorporated
already in our construction of $\{f,f\}_\ell$ as a comma object. 

Definition~\ref{defn:lax-morphism} can be modified to give other
flavours of weak morphism by replacing the operad $\{f,f\}_\ell$ by
analogous operads.

If we formed an {\em iso-comma object} $\{f,f\}\ps$ rather than the
comma object $\{f,f\}_\ell$ we would still obtain an operad, and
colax morphisms $\ct\to\{f,f\}\ps$ would then correspond to
pseudomorphisms of algebras, which are just the special case where
$\tilde{f}$ is invertible. 

If instead we formed the {\em pullback} $\{f,f\}_s$, then we would
obtain the strict morphisms, corresponding to the case where
$\tilde{f}$ is an identity. 

Alternatively, we could form the comma object $\{f,f\}_c$ in which the
direction of the 2-cell is reversed, and this would correspond to
reversing the direction of the $\tilde{f}$, and so to colax
morphisms. 

For each of these notions of morphism there is a corresponding notion
of 2-cell;  here we only  describe the case of 2-cells between lax morphisms. 
Suppose then that $\phi\colon f\to g$ is a natural transformation. As
well as $\{f,f\}_\ell$ we may form the operad $\{g,g\}_\ell$ and also
the object $\{f,g\}_\ell$ of $[\bbn,\Cat]$ appearing in the comma
object 
\[ \xymatrix{
\{f,g\}_\ell \ar[r]^c_{~}="1" \ar[d]_d^{~} & 
\<B,B\> \ar[d]^{\<f,B\>}_{~}\\
\<A,A\> \ar[r]_{\<A,g\>}^{~}="2" & \<A,B\> 
\ar@{=>}"1";"2" } \]
and now the pullback
\[ \xymatrix{
[\phi,\phi] \ar[r]^q \ar[d]_p & \{g,g\}_\ell \ar[d]^{\{\phi,g\}_\ell} \\
\{f,f\}_\ell \ar[r]_{\{f,\phi\}_\ell} & \{f,g\}_\ell } \]
and $[\phi,\phi]$ has a unique operad structure for which $p$ and $q$
are morphisms of operads. A  colax morphism $\ct\to[\phi,\phi]$
corresponds to a 2-cell between lax morphisms of  colax algebras. 

\begin{definition}
  We write $\clx\ct\Algl$ for the 2-category of colax
  \ct-algebras, lax morphisms, and algebra 2-cells, and we write
  $\nc\ct\Algl$ and 
  $\ct\Algl$ for the full sub-2-categories of normal colax algebras
  and strict algebras. Similarly
  we write $\clx\ct\Algs$, $\nc\ct\Algs$ and $\ct\Algs$ for the (non-full)
  sub-2-categories of strict morphisms. 
\end{definition}

There are analogous 2-categories with pseudomorphisms, designated with
a subscript ps, and of colax morphisms, with a subscript c. In many
cases there are parallel results for each such flavour of weak
morphism, and we write $\clx\ct\Algw$, $\nc\ct\Algw$, or $\ct\Algw$ if
we do not need to specify a particular choice. 

Before leaving this section we record the following standard fact. 

\begin{proposition}\label{prop:F*}
A morphism of operads $F\colon \ct\to\cq$ induces a 2-functor
$F^*\colon \cq\Algw\to\ct\Algw$. Similarly, 
a normal colax morphism $F\colon \ct\to\cq$ induces a 2-functor
$F^*\colon \nc\cq\Algw\to\nc\ct\Algw$, and a colax morphism $F\colon
\ct\to\cq$ induces a 2-functor $F^*\colon \clx\cq\Algw\to\clx\ct\Algw$.
\end{proposition}

\section{Background on the free skew monoidal category}

Skew monoidal structure on a category \cc involves certain basic functors
$\cc^n\to \cc$ and natural transformations between them, as well as
equations asserting that various functors and natural transformations,
derived from the basic ones via substitution, are equal. 

Thus there is an operad \cs whose (strict) algebras in \Cat are the
skew monoidal categories; more precisely, the corresponding 2-category
$\cs\Algs$ is isomorphic to the 2-category $\Skews$ of skew monoidal
categories, strict monoidal functors between them, and monoidal
transformations between these. 

Furthermore, it is not hard to check
that the lax morphisms of \cs-algebras are the (lax) monoidal functors
between skew monoidal categories, the colax morphisms are the
opmonoidal functors, and the pseudo morphisms are the strong monoidal
functors. Thus for each flavour w of weakness, there is an isomorphism
$\Skeww\cong \cs\Algw$.

As described in Section~\ref{sect:operads} above, the operad \cs can
equivalently be described as a club $\bs$. As an object of
$\Cat/\bbn$, this consists of the free \cs-algebra on the category
$1$, equipped with the unique strict morphism to $\bbn$ sending
the generator to $1\in\bbn$. 

In \cite{skewcoherence}, a specific construction of the free
\cs-algebra was given, under the name \Fsk.  As explained in the
introduction this construction was not fully explicit.  The goal of the
present section is to recall  the construction of
\cite{skewcoherence}, before 
giving a fully explicit description of \Fsk in Section~\ref{sect:revisited}.

%

\subsection{Ordinals} \label{sect:ordinals}

We write $\ord m$ for the ordinal $\{0,1,\ldots,m-1\}$. We can regard
$\ord m$ as a poset and hence as a category. 

A function
$\phi\colon \ord m\to \ord n$ is order-preserving if $i\le j$ implies
that $\phi(i)\le \phi(j)$. Thus the order-preserving functions between
ordinals are the functors. 

We can ask whether such functors have adjoints. Such a functor
$\phi\colon \ord m\to \ord n$ has a right adjoint if and only if it
preserves the least element: $\phi(0)=0$, as is always the case if
$\phi$ is surjective. When $\phi(0)=0$, the right
adjoint $\phi^*$ is given by 
\[ \phi^*(j) = \max\{ i \mid \phi(i)\le j\} \]
and may also be characterised by the fact that 
\[ \phi(\phi^*(j)) \le j < \phi(\phi^*(j)+1).\] 

In the context of ordinals, the usual adjointness property 
\[ \phi(i) \le j \Leftrightarrow i\le \phi^*(j) \]
can be expressed as 
\[ j < \phi(i) \Leftrightarrow \phi^*(j) < i. \]

It is useful to record:
\begin{proposition}
A right adjoint $\phi^*\colon \ord n\to \ord m$ itself has a right
adjoint $\phi_*$ if and only if $\phi(1)\neq 0$. In this case
\[ \phi_*(i) =
\begin{cases}
  \phi(i+1)-1 & \text{if $\phi(i)<\phi(i+1)$} \\
  \phi(i)-1 & \text{otherwise.}
\end{cases}
\]
\end{proposition}

\proof
We know that $\phi(0)=0$ since $\phi$ has a right adjoint. Now
$\phi^*$ will have a right adjoint if and only if $\phi^*(0)=0$; that
is, if $j>0$ implies that $\phi(j)>0$. But this will clearly follow
from the special case $\phi(1)>0$.

Now 
\[ \phi_*(i) = \max \{j\mid  \phi^*(j)\le i \}. \]
If $\phi^*(j)=i$ for some $j$, then the greatest such $j$ will clearly
be $\phi_*(i)$. We know that $\phi^*(j)=i$ if and only if $\phi(i)\le j<\phi(i+1)$,
which is possible if and only if $\phi(i)<\phi(i+1)$, and in that case
$\phi(i+1)-1$ will clearly be the greatest $j$. 

If $\phi(i)=\phi(i+1)$ then there is no $j$ with $\phi^*(j)=i$, so we
must settle for the greatest $j$ with $\phi^*(j)<i$, or equivalently
with $j<\phi(i)$. But the greatest such $j$ is clearly $\phi(i)-1$.
(Note that $\phi(0)<\phi(1)$, so $i=0$ is impossible, and so
$\phi(i)>0$ and $\phi(i)-1$ does exist. )
\endproof

\subsection{Left and right bracketing functions}\label{sect:lbfs}

The starting point for the description of the objects of \Fsk is the
Tamari lattice, which consists of all possible bracketings of an
$n$-fold product. These can be described explicitly using the idea of
a left bracketing function, given in \cite{TamariHuang}.

Let $\ord m=\{0,1,\ldots,m-1\}$ be a non-empty finite ordinal. A {\em left bracketing
  function}, or lbf, on $\ord m$ is a function $\ell\colon \ord m\to
\ord m$ satisfying three conditions:
\begin{enumerate}[(i)]
\item $\ell(j)\le j$ for all $j\in\ord m$
\item if $\ell(j)\le i<j$ then $\ell(j)\le \ell(i)$
\item $\ell$ preserves the top element $\top$ of $\ord m$.
\end{enumerate}
These are given the pointwise ordering: $\ell\le \ell'$ if and only if
$\ell(j)\le \ell'(j)$ for all $j$.

For example, in the case $\ord m=4$, this corresponds to the
bracketings of a 4-fold product, as in the diagram below 
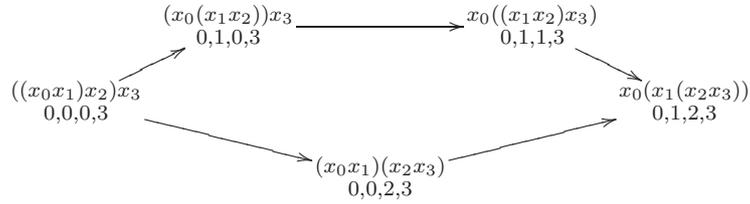
\[ \xy
(0,10)* { ((x_0x_1)x_2)x_3 \atop 0,0,0,3 }="a" ;
(20,20)* { (x_0(x_1x_2))x_3 \atop 0,1,0,3}="b" ; 
(60,20) * { x_0((x_1x_2)x_3) \atop 0,1,1,3}="c" ;
(80,10) * { x_0(x_1(x_2x_3)) \atop 0,1,2,3}="d" ;
(40,0) *  { (x_0x_1)(x_2x_3) \atop 0,0,2,3}="e" ;
\ar "a";"b"
\ar "b";"c"
\ar "c";"d"
\ar "a";"e"
\ar "e";"d"
\endxy \]
where, for example, the list $0,1,0,3$ denotes the lbf with
$\ell(0)=\ell(2)=0$, $\ell(1)=1$, and $\ell(3)=3$.

There are various ways to see the correspondence between bracketings
and lbfs.  Though not required in what follows let us give an example
illustrating one such way, which passes through the intermediate step
of a triangulation.  The bracketing $(x_0((x_1x_2)x_3))x_4$ of five elements
 corresponds to the triangulation of the $6$-gon as below.
$$
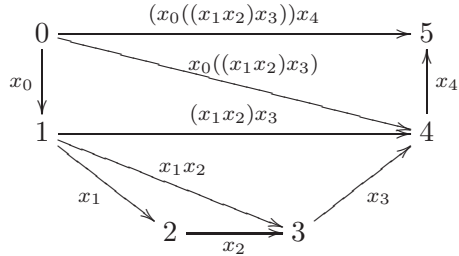
\xymatrix @C3pc {
0 \ar[d]_{x_0} \ar[drrr]^{x_0((x_1x_2)x_3)} \ar[rrr]^{(x_0((x_1x_2)x_3))x_4}  &&& 5 \\
1 \ar[dr]_{x_1} \ar[rrr]^{(x_1x_2)x_3} \ar[drr]^{x_1x_2} &&& 
4 \ar[u]_{x_4} \\
& 2 \ar[r]_{x_2} & 3 \ar[ur]_{x_3} 
}$$
The corresponding lbf $l\colon \ord 5 \to \ord 5$ is 
 obtained by defining $l(i)$ to be the least vertex of the triangle with middle $i+1$.  
 This only makes sense for $i \leq 3$; for the top element we are forced to define $l(4)=4$ so that
 the corresponding lbf is $0,1,1,0,4$.

We write $\Tam_m$ for the resulting poset. Its elements specify
bracketings of $m+1$-fold products. One advantage of the lbfs as a
description of the elements of $\Tam_m$ is that it makes it easy to
construct joins in $\Tam_m$. For if $\ell$ and $\ell'$ are lbfs then
so is the function $\ell\vee \ell'$ given by
$(\ell\vee\ell')(i)=\ell(i)\vee \ell'(i)$, where $\vee$ denotes the
join (maximum). This $\ell\vee\ell'$ is clearly the join of $\ell$ and
$\ell'$; thus $\Tam_m$ has binary joins. The function $\ord m\to\ord
m$ which has constant value 0 is an lbf, and is clearly the least
element of $\Tam_m$. Thus $\Tam_m$ has finite joins; but it is a
finite poset, so therefore has all joins and all meets. 

Every lbf $\ell\colon \ord m\to \ord m$ determines, and is determined
by, a function $r\colon \ord m\to \ord m$, connected via the relationships
\begin{align*}
   r(i) &= \min\{ j\mid \ell(j)<i\le j\} \\
\ell(i) &= \max\{i\mid i\le j<r(i)\} 
\end{align*}
Functions $r$ of this type are called {\em right bracketing
  functions} or rbfs. It turns out that if also $r'$ corresponds to
$\ell'$ then $\ell\le \ell'$ if and only if $r\le r'$. 

Since we often go back and forth between lbfs and rbfs, it is
convenient to introduce notation which is independent of this
choice. We therefore write $S$ for a particular element of the Tamari
poset, $\ell_S$ for the corresponding lbf and $r_S$ for the
corresponding rbf.

\subsection{Change of base for lbf's}
If $\sigma\colon\ord n \to \ord m$ is a surjective order preserving map
it has a right adjoint $\sigma^{*}$.  Proposition 5.6 of \cite{skewcoherence}
shows that if $l_{S}$ is an lbf on $\ord n$ the function $\sigma l_{S} \sigma^{*}\colon\ord m \to \ord m$ 
is an lbf on $\ord m$.  We write $\sigma S \sigma^{*}$ for the corresponding
element of $\Tam_{m}$.  

\subsection{The free skew monoidal category \Fsk on $1$}

We now turn to the construction of $\Fsk$ given in \cite{skewcoherence}.

\begin{definition}
  An object of \Fsk is a triple $(\ord m,u,S)$ where $\ord
m$ is a non-empty finite ordinal, $u$ is a subset of $\ord m$, and
$S\in\Tam_m$, with corresponding lbf $\ell_S\colon
\ord m\to \ord m$.
\end{definition}

This is thought of as an $m$-fold product, bracketed according to
$S$, with the generator $X$ in the positions specified by $u$,
and the unit $i$  elsewhere. 

The generating object is $(\ord 1,\ord 1,S)$ for the unique
$S\in\Tam_1$. For an arbitrary skew monoidal category \cc and object
$X\in\cc$, there is a unique strict monoidal functor $\Fsk\to\cc$
sending $(\ord 1,\ord 1,S)$ to $X$: this is what is meant by saying
that \Fsk is the free skew monoidal category on one object. 

To motivate the definition of the morphisms of \Fsk recall the
 category $\Ord_\bot$ of finite non-empty ordinals and functions preserving both order
and bottom element.  This admits a strictly associative skew monoidal
structure, with ordinal sum for tensor product and unit $1$.  

Accordingly the unit $\ord 1 \in \Ord_\bot$ determines a canonical map $\Fsk\to\Ord_\bot$ preserving
the skew monoidal structure strictly, and sending $(\ord m,u,S)$ to $\ord m$.
One of the main results of \cite{skewcoherence} is that this
$\Fsk\to\Ord_\bot$ is faithful; whereby morphisms $(\ord m,u,S) \to (\ord m,u,S)$ of $\Fsk$ can be identified
with certain morphisms $\ord m \to \ord m$ of $\Ord_\bot$.  

The question, then, is to identify which ones.  In \cite{skewcoherence} this was 
done in stages, starting with various special classes of morphism.  We begin with 
those corresponding to the associators $\alpha$.

\begin{definition}
If $(\ord n,u,S)$ and $(\ord n,u,T)$ are objects of $\Fsk$, we say that
the identity $1\colon\ord n \to \ord n$ defines a \emph{Tamari morphism} $1\colon(\ord n,u,S) \to (\ord n,u,T)$
if $S \leq T$.
\end{definition}

Next, those corresponding to applications of $\lambda$.

\begin{definition} 
A {\em shrink morphism} from $(\ord m,u,S)$ to $(\ord n,v,T)$ is an
order-preserving surjection $\sigma\colon \ord m\to \ord n$ satisfying
the following conditions:
\begin{enumerate}[(i)]
\item $\sigma$ and $\sigma^*$ restrict to mutually inverse bijections between
    $u$ and $v$
\item $\sigma\ell_S\sigma^*=\ell_T$
\item if $\sigma(j)=\sigma(j+1)$ then $\sigma(\ell_S(j))=\sigma(j)$.
\end{enumerate}
\end{definition}

\begin{remark}
As observed in \cite{skewcoherence}, if $\sigma^*$ restricts to a
bijection $v\to u$, the inverse is necessarily given by (the restriction
of) $\sigma$. But to say that $\sigma$ restricts to a bijection $u\to
v$ is not enough: we should also insist that $\sigma^*\sigma j=j$ for
any $j\in u$; in  other words, if $j\in u$ then $j$ is maximal
in the fibre $\sigma^{-1}\sigma(j)$.
\end{remark}

\begin{remark}\label{rmk:sigma-singular}
For an order-preserving surjection $\sigma$, to say that
$\sigma(j)=\sigma(j+1)$ is to say that $\sigma(j+1)\le \sigma(j)$, or
equivalently $j+1\le \sigma^*\sigma(j)$, or equivalently
$j<\sigma^*\sigma(j)$. Thus we can reformulate (iii) as 
\begin{enumerate}
\item[(iii)] if $j<\sigma^*\sigma(j)$ then $\sigma(\ell_S(j))=\sigma(j)$.
\end{enumerate}
\end{remark}

Combining the two classes yields a class of morphism named after 
the fact that they are precisely those morphisms of $\Fsk$ sent to surjections by the 
canonical map $\Fsk \to \Ord_\bot$.

\begin{definition}\label{def:surj}
An {\em $\Fsk$-surjection} $(\ord m,u,S)\to(\ord n,v,T)$ is an
order-preserving surjection $\sigma\colon \ord m\to \ord n$ that
factorizes as a Tamari morphism $1\colon(\ord m,u,S) \to (\ord m,u,S')$ followed by a shrink morphism 
$(\ord m,u,S') \to (\ord n,v,T)$.
\end{definition}

Corresponding to the application of $\rho$ there is the notion of a swell morphism.  
Combining these with the Tamari morphisms yields the $\Fsk$-injections, so named since they are precisely the maps sent to injections by the canonical $\Fsk \to {\Ord}_\bot$.

These are defined in \cite{skewcoherence} using duality, as below, but an elementary description
 can also be given -- see Section 10 of \cite{skewcoherence}.
For the definition using duality observe that given $(\ord n,v,T)$ of \Fsk we can
form the object $(\ord n\op,u,T\op)$ of \Fsk in which $T\op$ corresponds to the lbf
$r_T$ on $\ord n\op$.

\begin{definition}
An \Fsk-injection $(\ord n,v,T)\to(\ord m,u,S)$ is an order-preserving
left adjoint $\delta\colon \ord n\to \ord m$ for which $\delta^*\colon
\ord m\to \ord n$ defines an \Fsk-surjection $(\ord
m\op,u,S\op)\to(\ord n\op,v,T\op)$. Such a $\delta$ is a swell
morphism if and only if $\delta^*$ is a shrink morphism. 
\end{definition}

Finally we are in a position to describe the general case.

\begin{definition}
  A morphism in \Fsk from $(\ord m,u,S)$ to $(\ord n,v,T)$ is an
  order-preserving map $\phi\colon \ord m\to \ord n$ with a right
  adjoint that can be factorized as an $\Fsk$-surjection followed by
  an $\Fsk$-injection.
\end{definition}

\section{\Fsk revisited}\label{sect:revisited}

$\Fsk$-surjections, $\Fsk$-injections, and general morphisms of $\Fsk$ were defined
in terms of the existence of certain factorisations, which need not be unique.  In the 
present section we revisit each class, giving completely explicit descriptions of them.

\subsection{\Fsk surjections revisited}

In the definition of \Fsk-surjection the object $S^{\prime}$ is not uniquely determined.
The first step will be to describe a canonical choice for $S'$. It was
shown in \cite[Proposition~9.1]{skewcoherence} that there is a maximal choice
of $S'$; here we describe it more explicitly.

\begin{lemma}\label{lemma:9.1bis}
  Let $\sigma\colon \ord m\to \ord n$ be an order-preserving surjection and
  $\ell\colon \ord n\to \ord n$ an lbf. Consider the function
  $\ell^\sigma\colon\ord m\to \ord m$ given by 
\[ \ell^\sigma(j) = 
  \begin{cases}
    \sigma^*\ell \sigma(j) & \text{if $j=\sigma^*\sigma j$} \\
    j & \text{otherwise.} 
  \end{cases}
\]
 Then $\ell^\sigma$ is an
lbf and $\ell^\sigma\sigma^*=\sigma^*\ell$; thus also $\sigma\ell^\sigma\sigma^*=\ell$.
\end{lemma}

\proof
First observe that  if $j=\sigma^*h$, then $\sigma^*\sigma
j=\sigma^*\sigma\sigma^* h=\sigma^*h=j$, thus $j=\sigma^*\sigma j$ if and only if $j=\sigma^*h$
for some $h$. Thus
$\ell^\sigma\sigma^*h=\sigma^*\ell\sigma\sigma^*h=\sigma^*\ell h$ for
all $h$, and so $\ell^\sigma\sigma^*=\sigma^*\ell$ and therefore
$\sigma\ell^\sigma\sigma^*=\sigma\sigma^*\ell=\ell$. 

Thus we need only show that $\ell^\sigma$ is an lbf. 
If $j=\sigma^*\sigma j$ then $\ell^\sigma(j)=\sigma^*\ell\sigma(j)\le
\sigma^*\sigma(j)=j$, while otherwise $\ell^\sigma(j)=j$. Thus
$\ell^\sigma(j)\le j$ for all $j$.

Since $\sigma$, $\sigma^*$, and $\ell$ all preserve top elements, so
does $\ell^\sigma$.

Finally, suppose that $\ell^\sigma(j)\le i<j$. Then
$\ell^\sigma(j)\neq j$, so we must have $\sigma^*\sigma j=j$ and
$\ell^\sigma(j)=\sigma^*\ell\sigma j$.

If $i=\sigma^*\sigma i$ then $\ell^\sigma(j)=\sigma^*\ell\sigma j\le
\sigma^*\sigma i\le j=\sigma^*\sigma j$, and so applying $\sigma$
gives $\ell\sigma j\le \sigma i\le \sigma j$, and $\ell$ is an
lbf so $\ell\sigma j\le \ell\sigma i$, and 
\[ \ell^\sigma(j) = \sigma^*\ell\sigma(j)\le \sigma^*\ell\sigma(i)=\ell^\sigma(i).\]

If $i\neq \sigma^*\sigma i$ then $\ell^\sigma(i)=i$ and so $\ell^\sigma(j)\le
i=\ell^\sigma(i)$.
\endproof

If $T$ is the element of the Tamari lattice corresponding to $\ell$,
it is convenient to write $T^\sigma$ for the element of the Tamari lattice
corresponding to $\ell^\sigma$. 


\begin{proposition}\label{prop:Fsk-surj-fact1}
  Any \Fsk-surjection $\sigma\colon (\ord m,u,S)\to(\ord n,v,T)$
  factorises as 
\[ \xymatrix{
(\ord m,u,S) \ar[r]^-{1_{\ord m}} & (\ord m,u,S\vee T^\sigma) \ar[r]^-{\sigma}
& (\ord n,v,T) } \]
where the second factor is a shrink morphism. Furthermore, $S\vee
T^\sigma$ is the greatest $S'$ for which $\sigma\colon (\ord
m,u,S')\to (\ord n,v,T)$ is a shrink morphism. 
\end{proposition}

\proof
Since $S\le S\vee T^\sigma$, the first factor is an
\Fsk-surjection. We need to show that the second factor is a shrink
morphism and that $S\vee T^\sigma$ is maximal.

Since $\sigma\colon (\ord m,u,S)\to(\ord n,v,T)$ is an
\Fsk-surjection there is an $S'\ge S$ for which $\sigma\colon (\ord
m,u,S')\to(\ord n,v,T)$ is a shrink morphism. 

The fact that $\sigma^*$ restricts to a bijection $v\to u$ is
unchanged by passing from $S'$ to $\ell_S\vee T^\sigma$.

For the second condition we have
\begin{align*}
  \sigma\ell_{S\vee T^\sigma}\sigma^*j &= \sigma(\ell_S\sigma^*j\vee
  \ell_{T^\sigma}\sigma^*j) \\
  &= \sigma(\ell_S\sigma^*j\vee \sigma^*\ell_T j) \\
  &= \sigma\ell_S\sigma^*j\vee \sigma\sigma^*\ell_T j \\
  &= \sigma\ell_S\sigma^*j\vee \ell_T j \\
  &= \ell_T j
\end{align*}
since $\sigma\ell_S\sigma^*\le \sigma\ell_{S'}\sigma^*=\ell_T$.
Thus $\sigma\ell_{S\vee T^\sigma}\sigma^*=\ell_T$.

It remains to show that if $j<\sigma^*\sigma(j)$ then
$\sigma(\ell_{S\vee T^\sigma}(j))=\sigma(j)$. 
But if $j<\sigma^*\sigma(j)$ then $\ell^\sigma_T(j)=j$ and so $\ell_{S\vee
  T^\sigma}j=j$, and so $\sigma(\ell_{S\vee T^\sigma}j)=\sigma(j)$.

Now we show that $S\vee T^\sigma$ is the greatest $S'$ as in the
proposition. First observe that if $S_1$ and $S_2$ are any two such,
then $S_1\vee S_2$ is another, thus it will suffice to show that if
$S\vee T^\sigma<S'$ then $S'$ is not such an element. 

If $S\vee T^\sigma<S'$ then for some $j$ we have both
$\ell_S(j)<\ell_{S'}(j)$ and
$\ell_{T^\sigma}(j)<\ell_{S'}(j)$. Clearly this is impossible if
$\ell_{T^\sigma}(j)=j$, so we must have $j=\sigma^*\sigma(j)$ and
$\ell_{T^\sigma}(j)=\sigma^*\ell_T\sigma(j)$. Now
$\sigma^*\ell_T\sigma(j)<\ell_{S'}(j)$ and so 
\begin{align*}
  \ell_T(\sigma(j)) &< \sigma\ell_{S'}(j) \tag{adjointness} \\
  &= \sigma\ell_{S'}\sigma^*\sigma(j) \tag{$j=\sigma^*\sigma(j)$} \\
  &= \ell_T\sigma(j) \tag{shrink morphism} 
\end{align*}
giving a contradiction.  \endproof

We can now use this last result to provide a more explicit description
of $\Fsk$-surjections:

\begin{proposition}\label{prop:Fsk-surj}
  An order-preserving surjection $\sigma\colon \ord m\to\ord n$
  defines an \Fsk-surjection $(\ord m,u,S)\to(\ord n,v,T)$ if and only
  if
  \begin{enumerate}[(i)]
  \item $\sigma$ and $\sigma^*$ restrict to mutually inverse bijections between $u$ and $v$
    \item $\sigma\ell_S\sigma^*\le \ell_T$.
  \end{enumerate}
\end{proposition}

\proof
By the previous result, $\sigma$ will define an \Fsk-surjection if and
only if $\sigma\colon (\ord m,u,S\vee T^\sigma)\to (\ord n,v,T)$ is a shrink
morphism.

Condition (i) in the definition of shrink morphism is condition (i) in
the proposition. 
Condition (ii) in the definition of shrink morphism says that
$\sigma(\ell_S\vee \ell^\sigma_T)\sigma^*=\ell_T$. Now
\begin{align*}
  \sigma(\ell_S\vee \ell^\sigma_T)\sigma^* &= \sigma\ell_S\sigma^*
                                      \vee \sigma\ell^\sigma_T\sigma^* \\
           &= \sigma\ell_S\sigma^* \vee \sigma\sigma^*\ell_T  \\
  &= \sigma\ell_S\sigma^* \vee \ell_T 
\end{align*}
which is equal to $\ell_T $ if and only if condition (ii) in the
proposition holds.

Finally condition (iii) in the definition of shrink morphism says that
if $j<\sigma^*\sigma j$ then $\sigma(\ell_S\vee \ell^\sigma_T)j=\sigma
j$. But if $j<\sigma^*\sigma j$ then 
\begin{align*}
  \sigma(\ell_S\vee \ell^\sigma_T)j &= \sigma(\ell_S j \vee
                                      \ell^\sigma_T j) \\
  &= \sigma(\ell_S j\vee j) \\
  &= \sigma(j)
\end{align*}
and so this is automatic. 
\endproof

This in turn gives another factorisation:

\begin{proposition}\label{prop:Fsk-surj-fact2}
  Any \Fsk-surjection $\sigma\colon (\ord m,u,S)\to(\ord n,v,T)$
  factorises as 
\[ \xymatrix{
(\ord m,u,S) \ar[r]^-{\sigma} & (\ord n,v,\sigma S\sigma^*) \ar[r]^-{1}
& (\ord n,v,T). } \]
\end{proposition}

\proof 
The first factor satisfies the characterisation in
Proposition~\ref{prop:Fsk-surj}, so is an \Fsk-surjection. By that
same characterisation, $\sigma S\sigma^*\le T$, and so the second
factor is a Tamari morphism. 
\endproof

\subsection{\Fsk injections revisited}
We can deal quickly with $\Fsk$ injections using duality.  First we
dualise Lemma~\ref{lemma:9.1bis}.
%
%

\begin{lemma}
  Let $\delta\colon \ord n\to \ord m$ be a bottom-preserving injection
  with right adjoint $\delta^*$, and let $r\colon \ord n\to \ord n$ be
  an rbf. Consider the function $r^\delta\colon \ord m\to \ord m$
  given by 
\[ r^\delta(j) = 
  \begin{cases}
    \delta r\delta^*(j) & \text{if $j=\delta\delta^*(j)$} \\
    j & \text{otherwise.}
  \end{cases}
\]
Then $r^\delta$ is an rbf and $r^\delta \delta=\delta r$.
\end{lemma}

\proof
We can think of $r$ as an lbf on $\ord n\op$, and think of
$\delta^*$ as an order-preserving surjection $\ord m\op\to\ord n\op$
in which case $\delta$ becomes its right adjoint. Now apply Lemma~\ref{lemma:9.1bis}.
\endproof

Dualising the other results similarly, we have

\begin{proposition}\label{prop:Fsk-inj-fact1}
  Any \Fsk-injection $\delta\colon (\ord n,v,T)\to(\ord m,u,S)$
  factorises as 
\[ \xymatrix{
(\ord n,v,T) \ar[r]^-{\delta} & (\ord m,u,T^\delta\wedge S) \ar[r]^-{1}
& (\ord m,u,S) } \]
where the first factor is a swell morphism. (Furthermore
$T^\delta\wedge S$ is minimal with this property.)
\end{proposition}

\begin{proposition}\label{prop:Fsk-inj}
  An order- and bottom-preserving injection $\delta\colon
  \ord n\to\ord m$ defines an \Fsk-injection $(\ord n,v,T)\to (\ord
  m,u,S)$ if and only if 
  \begin{enumerate}[(i)]
  \item $\delta$ and $\delta^*$ restrict to mutually inverse bijections between $u$ and $v$
    \item $r_T\le \delta^* r_S\delta$ .
  \end{enumerate}
\end{proposition}

\begin{proposition}\label{prop:Fsk-inj-fact2}
  Any \Fsk-injection $\delta\colon (\ord n,v,T)\to(\ord m,u,S)$
  factorises as 
\[ \xymatrix{
(\ord n,v,T) \ar[r]^-{1} & (\ord n,v,\delta^*S\delta) \ar[r]^-{\delta}
& (\ord m,u,S). } \]
\end{proposition}

\subsection{General \Fsk morphisms revisited}  

In the definition of a general $\Fsk$-morphism, the underlying factorisation $\phi = \delta \circ \sigma$
in $\Ord_\bot$ must be the unique epi-mono factorisation.  So the definition
can be reformulated as follows.

\begin{definition}\label{defn:Fsk-general}
  A morphism in \Fsk from $(\ord m,u,S)$ to $(\ord n,v,T)$ is an
  order-preserving map $\phi\colon \ord m\to \ord n$ with a right
  adjoint, such that there exist an \Fsk-surjection 
\[ \xymatrix{ (\ord m,u,S) \ar[r]^-{\sigma} & (\im(\phi),\phi(u),R)
} \]
and an \Fsk-injection 
\[ \xymatrix{ (\im(\phi),\phi(u),R) \ar[r]^-{\delta} & (\ord n,v,T)
} \]
with $\phi=\delta\circ\sigma$ for some $R\in \Tam_{\im(\phi)}$.
\end{definition}

The $R$ appearing in the factorization need not be given explicitly.  In the following theorem we
show that there is a canonical choice for $R$, and use this to give the promised explicit
description of the morphisms of $\Fsk$.



\begin{theorem}\label{thm:Fsk}
  An order-preserving morphism $\phi\colon \ord m\to \ord n$ defines an \Fsk-morphism $(\ord
  m,u,S)\to (\ord n,v,T)$ if and only if
  \begin{enumerate}[(a)]
  \item $\phi$ has a right adjoint $\phi^*$
  \item $\phi$ and $\phi^*$ restrict to mutually inverse bijections between
    $u$ and $v$ 
    \item $\sigma\ell_S\sigma^*\le \delta^*\ell_T\delta_*$
  \end{enumerate}
where $\sigma\colon \ord m\to \im(\phi)$ and $\delta\colon
\im(\phi)\to\ord n$ are the induced maps. 
\end{theorem}

\proof
Given any $\phi\colon \ord m\to \ord n$ we may factorise it as a
surjection $\sigma\colon \ord m\to \im(\phi)$ followed by an injection
$\delta\colon \im(\phi) \to \ord n$. 

To say that $\sigma$ and $\sigma^*$ restrict to mutually inverse bijections
between $u$ and $\sigma(u)$ is to say that if $j\in u$ then
$\sigma^*\sigma(j)=j$, but
$\sigma^*\sigma(j)=\sigma^*\delta^*\delta\sigma(j)=\phi^*\phi(j)$, so
this says that if $j\in u$ then $\phi^*\phi(j)=j$.

Suppose that this is the case.  Then to say that $\delta$ and $\delta^*$ restrict to mutually inverse bijections
between $\sigma(u)$ and $v$ is then to say that $\delta$ maps
$\sigma(u)$ to $v$, and $\delta^*$ maps $v$ to $\sigma(u)$, and if
$i\in v$ then $\delta\delta^*(i)=i$.

Now $\delta$ maps $\sigma(u)$ to $v$ if and only if $\phi$ maps $u$ to
$v$. And $\delta^*$ maps $v$ to $\sigma(u)$ if and only if
$\sigma^*\delta^*$ maps $v$ to $\sigma^*\sigma(u)$; but
$\sigma^*\delta^*=\phi^*$ and $\sigma^*\sigma(u)=u$, thus this says
that $\phi^*$ maps $v$ to $u$. Also
$\delta\delta^*=\delta\sigma\sigma^*\delta^*=\phi\phi^*$. 

Thus to say that there are are mutually inverse bijections 
\[ \xymatrix{
u \ar@<1ex>[r]^{\sigma} & \phi(u) \ar@<1ex>[r]^{\delta}
\ar@<1ex>[l]^{\sigma^*} & v \ar@<1ex>[l]^{\delta^*} } \]
is just to say that condition (b) holds.

Suppose now that $R$ is given as in Definition~\ref{defn:Fsk-general}.
Since $\delta$ is an \Fsk-injection, we may use Proposition~\ref{prop:Fsk-inj-fact2}
to obtain a factorisation 
\[ \xymatrix{
(\ord m,u,S) \ar[r]^-{\sigma} & (\im(\phi),\phi(u),R)
\ar[r]^-{1} &  (\im(\phi),\phi(u),\delta^*T\delta)\ar[r]^-{\delta} &
(\ord n,v,T) } \]
where $r_{\delta^*T\delta}= \delta^* r_T\delta$. As observed in
\cite[Proposition~5.7]{skewcoherence}, the corresponding lbf
$\ell_{\delta^*T\delta}$ is $\delta^*\ell_T\delta_*$ where $\delta_*$
is the {\em right} adjoint of $\delta^*$. 
In this new factorisation, the middle factor is also an
\Fsk-surjection. Thus the composite of the first two factors is an \Fsk-surjection
and so 
\[ \xymatrix{
(\ord m,u,S) \ar[r]^-{\sigma} & (\im(\phi),\phi(u),\delta^*T\delta)\ar[r]^-{\delta} &
(\ord n,v,T) } \]
is also a factorization as in Definition~\ref{defn:Fsk-general}.

Thus we have proved that $\phi$ is a morphism in \Fsk if and only if,
in this last displayed composite, $\sigma$ is an \Fsk-surjection and
$\delta$ is an \Fsk-injection. By
Proposition~\ref{prop:Fsk-surj} this is equivalent to conditions (a),
(b), and (c).
\endproof 

By adjointness (c) above is equivalent to $\phi\ell_S\sigma^*\le \ell_T\delta_*$.  One should resist the temptation to use adjointness once again to
  transform the inequality $\phi\ell_S\sigma^*\le \ell_T\delta_*$ to
  $\phi\ell_S\le \ell_T\delta_*\sigma$. This would
  be valid if we knew that $\ell_S$ and $\ell_T$ were functors
  (order-preserving), but this need not be the case.
  
On the other hand, there is another possible reformulation:

\begin{theorem}\label{thm:Fsk2}
  An order-preserving morphism $\phi\colon \ord m\to \ord n$  defines an \Fsk-morphism $(\ord
  m,u,S)\to (\ord n,v,T)$ if and only if
  \begin{enumerate}[(a)]
  \item $\phi$ has a right adjoint $\phi^*$
  \item $\phi$ and $\phi^*$ restrict to mutually inverse bijections between
    $u$ and $v$ 
        \item[(d)] if $\phi(j)<\phi(j+1)$ then $\phi(\ell_S(j))\le \ell_T(\phi(j+1)-1)$.
  \end{enumerate}
\end{theorem}

\proof
Factorise $\phi$ as $\sigma\colon \ord m\to \im(\phi)$ and
$\delta\colon \im(\phi)\to \ord n$, as in Theorem~\ref{thm:Fsk}.
We need to show that (c) is equivalent to (d).

As mentioned above, (c) is equivalent to $\phi\ell_S\sigma^*\le \ell_T\delta_*$.  If $j=\sigma^*(h)$ then $h=\sigma\sigma^*(h)=\sigma(j)$. Hence $\phi\ell_S\sigma^*\le \ell_T\delta_*$ says
that if $j=\sigma^*(h)$ then $\phi\ell_S(j)\le
\ell_T\delta_*\sigma(j)$. Now $j=\sigma^*(h)$ for some $h$ if and only
if $j=\sigma^*\sigma(j)$; and this is equivalent, as in
Remark~\ref{rmk:sigma-singular} to $\sigma(j)<\sigma(j+1)$. This in
turn is clearly equivalent to $\phi(j)<\phi(j+1)$. 

Thus (c) is equivalent to the condition that if $\phi(j)<\phi(j+1)$
then $\phi\ell_S(j)\le \ell_T\delta_*\sigma(j)$. So we just need to
show that if $\phi(j)<\phi(j+1)$ then $\delta_*\sigma(j)=\phi(j+1)-1$.

Now $\delta_*\sigma(j)= \max\{ i \mid
\delta^*(i)\le\sigma(j)\}=\max\{i\mid \delta^*(i)=\sigma(j)\}$, since
$\delta^*$ is surjective and so there certainly exists an $i$ with
$\delta^*(i)=\sigma(j)$, and thus clearly the maximum must be of this
type. To say that $\delta^*(i)=\sigma(j)$ is to say that
\[ \delta(\sigma(j))\le i<\delta(\sigma(j)+1). \]

Since $\sigma(j)<\sigma(j+1)$ and $\sigma$ is surjective, we must have
$\sigma(j+1)=\sigma(j)+1$ thus the displayed inequality becomes 
\[ \phi(j) \le i < \phi(j+1) \]
and now the maximum value of $i$ is clearly $\phi(j+1)-1$.
\endproof

\section{Adjunctions of operads and skew monoidal categories as colax algebras}

In this section we describe adjunctions between the operad for skew monoidal
categories and other simpler operads $\ct$.  These adjunctions allow us to view skew monoidal
categories as colax $\ct$-algebras.  We begin by taking $\ct$ to be the terminal operad
before passing to another operad $\cl$ whose colax algebras fully capture skew monoidal structure.

\subsection{Colax and lax monoidal structure}

Let \cs be the operad for skew monoidal categories and \cn the operad
for strict monoidal categories. Since $\cn$ is the terminal operad, there
is a unique (strict) operad morphism
$P\colon \cs\to\cn.$ 
The induced 2-functors
$P^*\colon \cn\Algw\to\cs\Algw$
are the inclusions, for the various possible flavours of morphism, of strict monoidal categories in skew monoidal categories. 

In this section we will see that $P\colon \cs\to\cn$ has a \emph{colax} left adjoint, which
allows us to view each skew monoidal category as a colax $\cn$-algebra -- that is, a colax
monoidal category.

The universal property of the free skew monoidal category $\Fsk$ on $1$
gives a strict monoidal functor $\Fsk \to \bbn$, which sends $(\ord n,u,S)$ to
the cardinality $|u|$ of the subset $u$.  It follows that $\cs_m$ is the full subcategory
 of \Fsk consisting of all objects of the form $(\ord n,u,S)$ with $|u|=m$.  
 
\begin{theorem}\label{thm:P-left}
The map $P\colon \cs\to\cn$ has a left adjoint in $[\bbn,\Cat]$ with
identity unit.
\end{theorem}

\proof

This is equivalent to saying that $P\colon \cs_m\to\cn_m$ has a
left adjoint in \Cat with identity unit, for each $m\in\bbn$. Since
$\cn_m$ is the terminal category, this in turn is equivalent to
saying that each $\cs_m$ has an initial object. 


The initial object is $(\ord{m+1},\ord{m+1}/\{0\},\bot)$, where the specified
subset consists of all elements except $0$, and $\bot \in \Tam_{m+1}$ is the
bottom element of the Tamari lattice, corresponding to the lbf
$\ell_L$ with $\ell_L(m)=m$ and $\ell_L(i)=0$ if $i\neq m$.

We prove the universal property using the characterisation in 
Theorem~\ref{thm:Fsk2}. Suppose then that $(\ord n,v,T)\in\Fsk$ has
$|v|=m$. There is a unique order-preserving bijection $\theta\colon
\ord{m+1}/\{0\} \to v$. The only way to define a map $\phi\colon \ord{m+1}\to\ord{n}$
which preserves order and the bottom element, as required to have a right adjoint, and which restricts to $\theta$, is to define $\phi(i)=\theta(i)$ if $i\in u$, and $\phi(0)=0$. 

This proves uniqueness; it remains to verify the conditions (b) and
(d) in Theorem~\ref{thm:Fsk2}.  
 If $i \in \ord{m+1}/\{0\}$ then $\phi^{-1}\phi(i)=\{i\}$ 
so that $\phi^*$ restricts to give the inverse of $\theta$, as required for (b).
Finally we verify condition (d); that is, if $\phi(j)<\phi(j+1)$ then $\phi(\ell_L(j))\le
\ell_T(\phi(j+1)-1)$. 
But if $\phi(j)<\phi(j+1)$ then $j\neq m$ and so
$\phi(\ell_L(j))=\phi(0)=0\le \ell_T(\phi(j+1)-1)$ as required.
\endproof

Given a strict morphism $U\colon\ct \to \cq$ of \Cat-operads whose 
underlying morphism in $[\bbn,\Cat]$ has a left adjoint 
$F\colon\cq \to \ct$, the left adjoint $F$ admits the structure
of a colax morphism of operads: this is an instance of doctrinal adjunction \cite{Kelly-doctrinal}.
To describe the structure, let $\eta$ and $\epsilon$ denote the unit and
counit of the adjunction.  The components of the colax structure are given by
\begin{equation}
\xymatrix{
F(x \circ_i y) \ar[r]^-{{\eta_x} \circ_i {\eta_y}} & F(UFx \circ_i
UFy) = FU(Fx \circ_i Fy) \ar[r]^-{\epsilon_{x \circ_i y}} & Fx \circ_i Fy
}
\end{equation}
and
\begin{equation}
\xymatrix{
Fe_{\cq} \ar@{=}[r] & FUe_{\ct} \ar[r]^{\epsilon_{e_{\ct}}} & e_{\ct}  .
}
\end{equation}
where $e_{\ct} \in \ct_1$ and $e_{\cq} \in \cq_1$ are the units of the respective operads.

\begin{corollary}\label{cor:J}
  The left adjoint of Theorem~\ref{thm:P-left} is a colax morphism
  of operads, and so sends skew monoidal categories to colax monoidal
  categories. More precisely, it defines a 2-functor $J\colon\Skeww\to
  \clx\cn\Algw$ for each flavour $w$ of weak morphism. 
\end{corollary}
\proof
The colax structure follows as above. Composition with the left
adjoint therefore sends colax $\cs$-algebras to colax $\cn$-algebras,
and so in particular sends strict $\cs$-algebras to colax
$\cn$-algebras; that is, it sends skew monoidal categories to colax
monoidal categories.
\endproof

We may describe this process more explicitly. Let \cc be a skew
monoidal category.   This becomes colax monoidal when we define the
tensor product of the list $(a_1,\ldots,a_n)$  to be the tensor product
in \cc of $$ia_1\ldots a_n$$ bracketed to the left. Clearly this process loses structure: 
there is no way of recovering a general product $ab$  in \cc.  

In a moment we will describe another operad $\cl$, only slightly more
complex than $\cn$, whose colax algebras do encode the entire 
skew monoidal structure.  Before that, let us mention that there is a
 dual way of making a skew monoidal category into a lax
monoidal category. 

\begin{theorem}\label{thm:P-right}
The map $P\colon \cs\to\cn$ has a right adjoint in $[\bbn,\Cat]$ with
identity counit.
\end{theorem}

\proof
This amounts to proving that each $\cs_n$ has a terminal
object. Explicitly, this will be given by $(\ord{n+1},v,\top)$ where $v$
consists of all elements of $\ord{n+1}$ except the top, and $\top$ is the
greatest element of the Tamari lattice $\Tam_{n+1}$, with lbf $\ell_\top$
given by $\ell_\top(i)=i$ for all $i\in\ord{n+1}$.

But in fact there is no need to prove this separately; rather, we can
use the following duality argument. For any skew monoidal category
\cc, the opposite category $\cc\op$ is also skew monoidal when we use
the reverse tensor product: $a\ox_{\cc\op} b=b\ox_\cc a$; this also
interchanges the roles of $\lambda$ and $\rho$. This means that there
is an isomorphism $\cs\op_n\cong \cs_n$, and the image under this
isomorphism of the initial object of Theorem~\ref{thm:P-left} will be
terminal.
\endproof

The adjunction of  Theorem~\ref{thm:P-right} was established by Uustalu using a term rewriting approach in \cite{uustalu}.
By doctrinal adjunction we obtain:

\begin{corollary}
  The right adjoint of Theorem~\ref{thm:P-right} is a lax morphism
  of operads, and so sends skew monoidal categories to lax monoidal
  categories. \endproof
\end{corollary}

This time the product in the lax monoidal category of the list $(a_1,\ldots,a_n)$ 
is the tensor product $$a_1\ldots a_ni$$ in the skew monoidal category,
bracketed to the right. 

\subsection{Colax \cl-algebras}

As we have mentioned, the passage from a skew monoidal category to the 
associated (co)lax monoidal category loses information.  In order to rectify
this problem, we may consider intermediate structures between strict monoidal
and skew monoidal categories, in the following sense.  Suppose that \cl is an operad, and that $P\colon \cs\to\cn$ factorises as 
\[ \xymatrix{
\cs \ar[r]^Q & \cl \ar[r]^R & \cn } \]
so that $P^*\colon \cn\Algw\to\cs\Algw$ factorises as 
\[ \xymatrix{
\cn\Algw \ar[r]^{R^*} & \cl\Algw \ar[r]^{Q^*} & \cs\Algw } \]
where w could be any of s, $\ell$, c, or ps. 

If each $Q\colon \cs_m\to\cl_m$ has a left adjoint $F$, then the
various $F$ inherit the structure of a colax morphism $\cl\to\cs$ of
operads, and so composition with $F$ sends skew monoidal categories to
colax \cl-algebras. We shall apply this for a specific choice of \cl, whose algebras
will be the following structures.



\begin{definition}
  A $\lambda$-algebra is a skew monoidal category for which both the
associativity maps $\alpha$ and the right unit maps $\rho$ are
identities.
\end{definition}

This can be considered as a structure in its own right:
it is a category \cc equipped with a strictly associative
multiplication $\cc\x\cc\to\cc$ and a strict right unit $i$; there is
also a natural transformation $\lambda\colon ia\to a$ satisfying three
conditions: $\lambda_a\ox b=\lambda_{a\ox b}$, $a\ox \lambda_b=1$, and
$\lambda_i=1$. As such it is clear that $\lambda$-algebras are the algebras
for an operad which will be called \cl. This operad can be described as follows.


\begin{proposition}
The operad \cl for $\lambda$-algebras has $\cl_0=\{l\}$ whilst $\cl_n$ is the two-element poset $\{\myl \le
\myt \}$ for $n>0$.
The multiplication $\cl_n\x \cl_{k_1}\x\ldots\x \cl_{k_n}\to \cl_{k_1+\ldots+k_n}$ is given by 
 \[ x(x_1\ldots,x_n) =
 \begin{cases}
   \myt & \text{if $x,x_1=\myt$} \\
   \myl & \text{otherwise}
 \end{cases}
\]
and the unit by $\myt\in \cl_1$.
\end{proposition}

In what follows we will often write $\myl_n$,$\myt_n$ to indicate that we are referring to $\myl,\myt \in \cl_n$.

\begin{proof}
Let $\cc$ be a category and consider $L\cc=\sum_{n\in\bbn} \cl_n
\times \cc^n$.  We write $\overline{a}=(a_1,\ldots,a_n)$ for a
typical element of $\cc^n$, and $\overline{a}\overline{b}$ for the
concatenation of lists $\overline{a}$ and $\overline{b}$. We equip
$L\cc$ with multiplication $(x_m,\overline{a})
\otimes (y_n,\overline{b}) = (x_{m+n},\overline{a},\overline{b})$ and
unit $(l_{0},-)$.  The left unit maps are the morphisms $
(\myl_0,-)\otimes (x_n,\overline{a}) = (\myl_n,\overline{a}) \to
(x_{n},\overline{a})$ induced by $l_n \leq x_n$.

Next we show that $L\cc$ is the free $\lambda$-algebra on $\cc$.  To
this end, consider a $\lambda$-algebra $\cd$ and functor $F\colon\cc \to
\cd$.  We must show that there is a unique structure-preserving
morphism $ \widehat{F}   \colon L\cc \to\cd$ sending each $(\myt_1,a)$ to
$Fa$.  This is straightforward. We can and must define
$ \widehat{F}  (\myt_n,\overline{a})$ as the $n$-fold tensor product
$\otimes_{i=1}^{n}Fa_i$, with $ \widehat{F}  (\myl_n,\overline{a})= i  \otimes
(\otimes_{i=1}^{n}Fa_i)$; while $ \widehat{F}  $ applied to a morphism $(\leq,\overline{a})\colon (\myl_n,\overline{a}) \to  (\myt_n,\overline{a})$ is the left unit map $ i \otimes (\otimes_{i=1}^{n}Fa_i) \to (\otimes_{i=1}^{n}Fa_i)$.

The unique $\lambda$-algebra map $L1 \to \bbn$ sending the generator
to $1 \in \bbn$ produces the values of our operad $\cl_i$ as its
fibres.  The components of the multiplication are calculated as the
components of the counit $1^*_{L1}\colon LL1 \to L1$.
\end{proof}

There is a unique (strict) operad morphism $R\colon \cl\to\cn$ sending
$\xi_n$ to $n$,, and the induced 
2-functor $R^*\colon \cn\Algw\to\cl\Algw$ is
the inclusion of strict monoidal categories in $\lambda$-algebras. 
Though we will not use this fact, we note that $R\colon \cl\to\cn$ has both adjoints in $[\bbn,\Cat]$ since each
$\cl_n$  has both an initial and a terminal object.

Of more interest to us is the unique (strict) operad morphism $Q\colon \cs\to\cl$ for which the
induced 2-functor $Q^*\colon \cl\Algw\to\cs\Algw$ is the inclusion of
\cl-algebras in skew monoidal categories. Explicitly, $Q\colon
\cs_n\to\cl_n$ sends $(\ord n,u,S)$ to $\myt_{|u|}$ if the bottom
element of $\ord n$ is in $u$, and $\myl_{|u|}$ otherwise.

\begin{theorem}\label{thm:Q-left}
The map $Q\colon \cs\to\cl$ has a left adjoint $H$ in $[\bbn,\Cat]$ with
identity unit.
\end{theorem}

\proof
We need to show that each $Q\colon\cs_m\to\cl_m$ has a left adjoint
with identity unit. Since $\myl_m \in \cl_m$ is initial we define $H\myl_m$ to 
be the initial object $(\ord m,\ord{m+1}\setminus\{0\},\bot)$ -- constructed
in Theorem~\ref{thm:P-left}.  Then $H\myl_m$ has the correct universal property.
By construction $QH\myl_m=\myl_m$ so that the unit component is the identity at $\myl_m$.

That leaves the case $\xi=\myt_m$. In this case we shall show that
$H\myt_m=(\ord m,\ord m,\bot)$, where $\bot\in\Tam_m$ is the bottom
element. 
If $(\ord n,v,T)\in\cs_m$ then $|v|=m$, and $\myt \le Q(\ord n,v,T)$ just when $0\in v$.  
There is then a unique order-preserving bijection $\theta\colon\ord m\to v$, and composing this with the inclusion
$v\to\ord n$ gives the order-preserving $\phi\colon \ord m\to\ord n$
mapping $\ord m$ bijectively to $v$; furthermore, since $0\in v$ it
satisfies $\phi(0)=0$ and so has a right adjoint $\phi^*$.  
This restricts to $\theta^{-1}$ since the given subset of $\ord m$ is its entirety.
Finally it satisfies
$\phi(\ell_L(j))\le \ell_T(\phi(j+1)-1)$ for all $j$ not equal to the
top element of $\ord m$, and so defines a morphism in \Fsk.
The unit at $\myt_m$ is once again the identity.
\endproof

\begin{remark}\label{remark:injection}

There is a little more we can say about morphisms $Hx_m \to (\ord
n,v,T) \in \cs_m$ in the context of the 
adjunction $H \dashv Q$.  First, such a morphism is unique if it exists, since $\cl_m$ is a poset.
Second, tracing through the construction of the adjunction we see that
the unique map $Hx_m \to (\ord n,v,T)$ in $\cs_m$ corresponding to an identity $x_m=Q(\ord n,v,T)$
is an $\Fsk$-injection.  In other words, the components of the counit are $\Fsk$-injections.
%
%
\end{remark}

\begin{corollary}\label{cor:H}
The left adjoint $H$ of Theorem~\ref{thm:Q-left} is a normal colax morphism of
operads, and so sends skew monoidal categories to normal colax
\cl-algebras. More precisely, it defines a 2-functor
$H^*\colon\Skeww\to\nc\cl\Algw$ for any flavour $w$ of weak morphism. 
\end{corollary}

\proof We just need to check that the colax morphism $H$ is
normal. $H$ sends the unit $t_1\in \cl_1$ of $\cl$ to the unit
$(\ord 1,\ord 1,\bot)$ of $\cs$, and this object has no non-identity
endomorphisms, so the colax structure map $H(t_1)\to (\ord 1,\ord 1,\bot)$ can
only be the identity.
\endproof

We shall see in Section~\ref{sec:skewmonoidal} that this 2-functor is fully faithful, and we shall
also characterise its image. 

\begin{remark}
Between the operads \cs and \cn there are various other possible operads one may consider.  Though not necessary in what follows let us briefly mention a fuller picture of such possibilities.
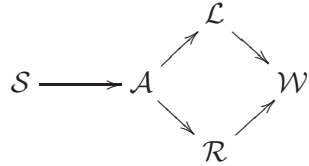
\begin{equation*}\label{eq:operadmaps}
\xymatrix @R1pc @C1pc { &&& \cl  \ar[dr] \\
\cs  \ar[rr] && \ca \ar[ur]  \ar[dr] && \cn \\
&&& \crr \ar[ur]
}
\end{equation*}
Here, in addition to \cs, \cl and \cn, are the operads \ca for skew
monoidal categories in which $\alpha$ is an identity and \crr for skew
monoidal categories in which both $\alpha$ and $\lambda$ are
identities. 
(In fact $\crr$ is dual to $\cl$, in the sense that $\crr_n =
\cl_n^{op}$.)  
In this diagram all of the morphisms on or above the horizontal have left adjoints in $[\bbn,\Cat]$ with identity unit whilst all those on or below the horizontal have right adjoints in $[\bbn,\Cat]$ with identity counit.
\end{remark}

\section{\LNC-algebras and \LNC-morphisms}

We have seen that each skew monoidal category gives rise to a normal colax $\cl$-algebra.
In this section we identify the property that characterises the colax
$\cl$-algebras arising in this way: we call such objects
\LNC-algebras; the ``LB'' stands for left-bracketed and the ``C'' for
colax.

We give a detailed analysis of the corresponding notion of \LNC-morphism of operads, which we will use in the following section to establish the correspondence
with skew monoidal categories.

Recall that $\cl_n=\{\myt_n,\myl_n\}$ for each $n\ge 0$ and
$\cl_0=\{\myl_0\}$.  The multiplication for \cl satisfies $\myt_2\circ_1x_n=x_{n+1}$ for all $x_n\in \cl_n$. 

Consider a normal colax \cl-algebra. This consists of a category $A$
equipped with functors $\myl=\myl_n\colon A^n\to A$ for each $n\ge 0$ and
$\myt=\myt_n\colon A^n\to A$ for each $n>0$, a natural transformation
$\lambda=\lambda_n\colon \myl_n\to\myt_n$ for each $n>0$, and 
suitably coherent natural transformations $\Gamma_{x,i,y} \colon m_{x\circ_i y}\to m_x
\circ_i m_y$ for each $x\in \cl_k$, $y\in \cl_n$, and
$i\in\{1,\ldots,k\}$.

\begin{definition}
  We say that a normal colax \cl-algebra $A$ is an {\em \LNC-algebra} if each $\Gamma_{\myt_2,1,x}\colon m_x\to
  m_{\myt_2}\circ_1 m_x$ is an identity.
\end{definition}

In particular, there are equalities 
\[ m_x(a_1,\ldots,a_{n+1})=m_{\myt_2}(m_x(a_1,\ldots,a_n),a_{n+1}). \]
We write $\LNC\Algw$ for the full sub-2-category of $\nc\cl\Algw$
consisting of the \LNC-algebras.

Since normal colax algebra structure on $A$ amounts to a normal colax
morphism $\ct\to\End(A)$, there is a natural extension of the previous
definition. 

\begin{definition}
  For an operad \ct, an {\em \LNC-morphism} from \cl to \ct is a normal colax
  morphism $F\colon \cl\to\ct$ for which each
  \[  \tilde{F}\colon  F(\myt_2\circ_1 x_n)\to F(\myt_2)\circ_1 F(x_n) \]
  is an
  identity. 
\end{definition}

Thus an \LNC-morphism $\cl\to\End(A)$ is the same as an \LNC-algebra
structure on $A$.

\begin{example}\label{ex:H}
The normal colax morphism $H\colon \cl\to\cs$ of Corollary~\ref{cor:H}
is an \LNC-morphism.  
First observe that 
\begin{align*}
  H(\myt_2)\circ_1 H(\myl_n) &= (\ord 2,\ord 2,\bot_2)\circ_1 (\ord n,
                               \ord n\setminus \{0\}, \bot_n) \\
  &= (\ord{n+1},\ord{n+1}\setminus\{0\},\bot_{n+1}) \\
  &= H(\myl_{n+1}) \\
  &= H(\myt_2\circ_1 \myl_n).
\end{align*}
The colax structure map $H(\myt_2\circ_1\myl_n)\to H(\myt_2)\circ_1
H(\myl_n)$ is given by the composite 
\[ \xymatrix{
H(\myt_2\circ_1 \myl_n) = H(QH\myt_2\circ QH\myl_n) =
HQ(H\myt_2\circ_1 H\myl_n) \ar[r] & H\myt_2\circ_1 H\myl_n } \]
where the first equality is because the unit of $H\dashv Q$ is an
identity, and the second because $Q$ is a (strict) morphism of
operads, while the unnamed arrow is the component at $H\myt_2\circ
H\myl_n$ of the counit. But this is also the component of the counit
at $H(\myl_{n+1})$, which is an identity by one of the triangle
equations. The proof that the colax structure map
$H(\myt_2\circ_1\myt_n)\to H(\myt_2)\circ_1 H(\myt_n)$ is an identity
is similar. 
\end{example}

Note that unlike colax morphisms, \LNC-morphisms are only defined when
the domain is \cl. However, if $F\colon \cl\to\ct$ is an
\LNC-morphism and $G\colon \ct\to\cq$ is a (strict) morphism of
operads, then the composite $G\circ F$ is also an \LNC-morphism. 

We now analyse what exactly is involved in giving an \LNC-morphism.

\begin{proposition}
  For an \LNC-morphism $F\colon \cl\to\ct$, the maps
  \[ \tilde{F}_{\myt_n,1,x_k}\colon  F(x_{n+k-1}) \to F(t_n)\circ_1 F(x_k)\]
  are also identities for all $n$.
\end{proposition}

\proof
The case $n=1$ holds by counitality, and the case $n=2$ by the \LNC\ 
condition. For $n>2$, use the coassociativity condition
\[
  \xymatrix{ 
  F(x_{n+k-1}) \ar[r]^{\tilde{F}_{t,1,x}}
  \ar@{=}[d]_{\tilde{F}_{t,1,x}} & F(t_n)\circ_1 F(x_k)
  \ar@{=}[d]^{\tilde{F}_{t,1,t}\circ_1 1} \\
  F(t_2)\circ_1 F(x_{n+k-2}) \ar[r]_{1\circ_1 \tilde{F}_{t,1,x}}
  & F(t_2)\circ_1 F(t_{n-1})\circ_1 F(x_k) } \]
and induction.
\endproof

\begin{proposition}
In an \LNC-morphism $F\colon\cl\to \ct$,  all of the functors $\cl_n\to \ct_n$ are
determined by
$F(\myt_2)\in \ct_2$, $F(\myl_0)\in \ct_0$, and $F(\mylam_1)\colon
F(\myl_1)\to F(\myt_1) \in \ct_1$.
\end{proposition}

\proof
This follows by a straightforward induction using the fact that
$\myt_{n+1}=\myt_2\circ_1\myt_n$,
$\myl_{n+1}=\myt_2\circ_1\myl_n$, and
$\mylam_{n+1}=\myt_2\circ_1\mylam_n$.
\endproof

\begin{lemma}\label{lemma:fullP}
For an \LNC-morphism $F\colon \cl\to \ct$, all of the colax structure
maps 
$\tilde{F}_{x,j,y}\colon F(x\circ_j y)\to F(x)\circ_j
F(y)$ are determined by those for which $x=\myt\in \cl_n$ and $j=n$.
\end{lemma}

\proof
Let $x\in \cl_n$ and $y\in \cl_k$ be given and consider
$\tilde{F}_{x,j,y}\colon F(x\circ_j y)\to F(x)\circ_j F(y)$. If $j=1$ and $x=\myt$ then this is an
identity. 

\subsubsection*{Step 1: $j>1$.} In this case $x\circ_j y=x$, regardless
of the value of $y$. By coassociativity, the diagram 
\[ \xymatrix {
F(x_{n+k-1}) \ar[rr]^{\tilde{F}_{x_n,j,y_k}}
  \ar@{=}[dd]_{\tilde{F}_{(\myt_{j+k-1},1,x_{n-j+1})}} &&
F(x_n)\circ_j F(y_k)
\ar@{=}[d]^{\tilde{F}_{(\myt_{j},1,x_{n-j+1})}\circ_j 1} \\
&& \left( F(\myt_{j})\circ_{1} F(x_{n-j+1}) \right) \circ_j F(y_k)
\ar@{=}[d]  \\
F(\myt_{j+k-1})\circ_1 F(x_{n-j+1})
\ar[rr]_-{\tilde{F}_{\myt_j,j,y_k}\circ 1}
&& \left( F(\myt_{j})\circ_j F(y_{k}) \right)\circ_1 F(x_{n-j+1}) } \]
commutes and by the \LNC\ property the verticals are identities. Thus the upper
horizontal is equal to the lower horizontal, which depends only on the
$\tilde{F}$ of the given form. 

\subsubsection*{Step 2: $j=1<n$.}
If $x=\myt$ there is nothing to prove, so suppose that $x=\myl$. 
In this case $x_n\circ_j y_k=\myl_{n+k-1}$ regardless of the value of
$y$. If $n>1$ then by coassociativity, the diagram
\[ \xymatrix @C5pc {
F(\myl_{n+k-1}) \ar@{=}[r]^{\tilde{F}_{\myt_2,1,\myl_{n+k-2}}}
\ar[d]_{\tilde{F}_{\myl_n,1,y_k}} &
F(\myt_2)\circ_1 F(\myl_{n+k-2}) \ar[d]^{1\circ_1\tilde{F}_{\myl_{n-1},1,y_k}}
\\
F(\myl_n)\circ_1 F(y_k)
\ar@{=}[r]_-{\tilde{F}_{\myt_2,1,\myl_{n-1}}\circ_1 1} & 
F(\myt_2)\circ_1 F(\myl_{n-1})\circ_1 F(y_k) } 
 \] 
commutes. We may now use induction to reduce to the case
$n=1$. In that case use coassociativity as in
\[ \xymatrix{
F(\myl_k) \ar[rr]^-{\tilde{F}_{\myl,1,y}}
\ar@{=}[dd]_{\tilde{F}_{\myt_{k+1},1,\myl_0}} &&
F(\myl_1)\circ_1 F(y_k) \ar@{=}[d]^{\tilde{F}_{\myt_2,1,\myl_0}\circ_1
1} \\
&& \left( F(\myt_2)\circ_1 F(\myl_0)\right) \circ_1 F(y_k) \ar@{=}[d]
\\
F(\myt_{k+1})\circ_1 F(\myl_0) \ar[rr]_-{\tilde{F}_{\myt_2,2,y_k}\circ_1
1} &&
\left( F(\myt_2)\circ_2 F(y_k) \right) \circ_1 F(\myl_0) 
    } \]
to reduce to $\tilde{F}_{\myt_2,2,y_k}$.
  \endproof

\begin{proposition}\label{prop:Gamma-reduction}
The $\tilde{F}$'s are determined by $\tilde{F}_{\myt_2,2,\myt_2}$ and $\tilde{F}_{\myt_2,2,\myl_0}$.
\end{proposition}

\proof
We have already reduced to the case of $\tilde{F}_{\myt_n,n,y_k}$;
and the case $n=1$ is already covered by
counitality. If $n>2$ use coassociativity as in 
\[ \xymatrix @C4pc {
    F(\myt_{n+k-1}) \ar[r]^{\tilde{F}_{\myt_n,n,y_k}} \ar@{=}[dd]_{\tilde{F}_{\myt_{k+1},1,\myt_{n-1}}} &
    F(\myt_n)\circ_n F(y_k) \ar@{=}[d]^{\tilde{F}_{\myt_{2},1,\myt_{n-1}}\circ_n 1} \\
      & \left( F(\myt_{2})\circ_{1} F(\myt_{n-1})\right) \circ_n F(y_k)
      \ar@{=}[d] \\
      F(\myt_{k+1})\circ_{1} F(\myt_{n-1})
      \ar[r]_-{\tilde{F}_{\myt_2,2,y_k}\circ_1 1} & 
        \left( F(\myt_{2})\circ_{2} F(y_k) \right)  \circ_1 F(\myt_{n-1}) }
    \]
    to reduce to the case where $n=2$. 

If now $y_k\in \{\myl_0,\myt_2\}$ there is nothing to prove; if
$y_k=\myt_1$ use counitality; otherwise use coassociativity as in 
\[ \xymatrix{
F(\myt_{k+1}) \ar[rr]^{\tilde{F}_{\myt_2,2,y_{k}}}
\ar[dd]_{\tilde{F}_{\myt_3,2,y_{k-1}}} && 
F(\myt_2)\circ_2 F(y_{k}) \ar@{=}[d]^{1 \circ_2 \tilde{F}_{\myt_2,1,y_{k-1}}} \\
&& F(\myt_2)\circ_2 \left(F(\myt_2)\circ_1 F(y_{k-1}) \right) \ar@{=}[d] \\
F(\myt_3)\circ_2 F(y_{k-1})
\ar[rr]_-{\tilde{F}_{\myt_2,2,\myt_2}\circ_2 1 } &&
(F(\myt_2)\circ_2 F(\myt_2) )\circ_2 F(y_{k-1}) 
} \]
to reduce to $\tilde{F}_{\myt_2,2,\myt_2}$ and
$\tilde{F}_{\myt_3,2,y_{k-1}}$, and now repeat the process in Step 1 of
Lemma~\ref{lemma:fullP} to deal with deal with $\tilde{F}_{\myt_3,2,y_{k-1}}$.
\endproof

\begin{proposition}\label{prop:spec-gen}
  All of the structure of an \LNC-morphism $F\colon \cl\to \ct$ is
  determined by $F(\myt_2)$, $F(\myl_0)$, $F(\mylam_1)$,
  $\tilde{F}_{\myt_2,2,\myt_2}$, and $\tilde{F}_{\myt_2,2,\myl_0}$.
\end{proposition}

\proof
This follows directly from the previous results.
\endproof

\section{Skew monoidal categories as \LNC-algebras}\label{sec:skewmonoidal}

In this section we describe the perfect correspondence between skew monoidal categories
and \LNC-algebras.

\subsection{From an \LNC-algebra to a skew monoidal category}
In Proposition~\ref{prop:spec-gen}, we saw that a special morphism
$F\colon \cl\to \ct$ is determined by a small amount of data. Next we
apply this to the case where $T$ is $\End(A)$, or $\{F,F\}\wk$, or $[\rho,\rho]$.

Suppose that $A$ is a normal colax $\cl$-algebra satisfying Property
\LNC. The corresponding \LNC-morphism $F\colon \cl\to \End(A)$ is
determined by:
\begin{itemize}
\item $m=F(\myt_2)\colon A^2\to A$
\item $i=F(\myl_0)\in A$
\item $\lambda=F(\mylam_1)\colon m\circ_1 i\to 1$
\item $\alpha=\tilde{F}_{\myt_2,2,\myt_2}\colon m\circ_1 m\to m\circ_2
  m$
\item $\rho=\tilde{F}_{\myt_2,2,\myl_0}\colon 1\to m\circ_2 i$.
\end{itemize}
We shall show that $(m,i,\alpha,\lambda,\rho)$ satisfy the five axioms
needed to define a skew monoidal structure on $A$. 

\begin{example}\label{ex:a.r1}
Consider $\tilde{F}_{\myt_2,2,\myl_1}$. Observe that
$\myl_1=\myt_2\circ_1\myl_0$. By coassociativity the diagram 
\[ \xymatrix{
F(\myt_2) \ar[rr]^{\tilde{F}_{\myt_2,2,\myl_1}}
\ar[d]_{\tilde{F}_{\myt_3,2,\myl_0}} &&
F(\myt_2)\circ_2 F(\myl_1)
\ar@{=}[d]^{1\circ_2\tilde{F}_{\myt_2,1,\myl_0}} \\
F(\myt_3)\circ_2 F(\myl_0)
\ar[r]_-{\tilde{F}_{\myt_2,2,\myt_2}\circ_2 1} & \left(F(\myt_2)\circ_2F(\myt_2)\right)\circ_2F(\myl_0)  \ar@{=}[r] & F(\myt_2)\circ_2(F(\myt_2)\circ_1F(\myl_0))
} \]
commutes. The lower horizontal is $\alpha\colon (xi)y\to x(iy)$. The left
vertical is $m\circ_1 \tilde{F}_{\myt_2,2,\myl_0}$, which is
$\rho1\colon xy\to (xi)y$. 
\end{example}

\begin{example}\label{ex:a.a1}
  Consider $\tilde{F}_{\myt_3,2,\myt_2}$. By coassociativity the
  diagram 
\[ \xymatrix{
F(\myt_{4}) \ar[rr]^{\tilde{F}_{\myt_2,2,\myt_{3}}}
\ar[d]_{\tilde{F}_{\myt_3,2,\myt_{2}}} && 
F(\myt_2)\circ_2 F(\myt_3) \ar@{=}[d]^{1 \circ_2
  \tilde{F}_{\myt_2,1,\myt_2}} \\
F(\myt_3)\circ_2 F(\myt_{2})
\ar[r]_-{\tilde{F}_{\myt_2,2,\myt_2}\circ_2 1} &
(F(\myt_2)\circ_2 F(\myt_2) )\circ_2 F(\myt_2) \ar@{=}[r] &
F(\myt_2)\circ_2(F(\myt_2)\circ_1 F(\myt_2)) } \]
commutes. The lower horizontal has the form $a\colon (w(xy))z\to
w((xy)z)$, and the left vertical is $a1\colon ((wx)y)z\to (w(xy))z$. 
\end{example}

\begin{proposition}\label{prop:ar}
The diagram 
\[ \xymatrix{
xy \ar[r]^{\rho} \ar[dr]_{1\rho} & (xy)i \ar[d]^{\alpha} \\ & x(yi) } \]
commutes.
\end{proposition}

\proof
By coassociativity the diagram 
\[ \xymatrix{
F(\myt_2) \ar[rr]^{\tilde{F}_{\myt_3,3,\myl_0}}
\ar@{=}[d]_{\tilde{F}_{\myt_2,2,\myt_1}} && 
F(\myt_3)\circ_3F(\myl_0)
\ar[d]^{\tilde{F}_{\myt_2,2,\myt_2}\circ_3 1} \\
F(\myt_2)\circ_2F(\myt_1)  \ar[r]_-{1  \circ_2\tilde{F}_{\myt_2,2,\myl_0}} &
F(\myt_2)\circ_2(F(\myt_2)\circ_2F(\myl_0))
\ar@{=}[r] & (F(\myt_2)\circ_2F(\myt_2))\circ_3F(\myl_0) } \]
commutes. This agrees with the diagram in the proposition (use
Proposition~\ref{prop:Gamma-reduction} to identify the top row with
$\rho \colon xy\to (xy)i$.)
\endproof

\begin{proposition}\label{prop:al}
The diagram 
\[ \xymatrix{
(ix)y \ar[r]^{\alpha} \ar[dr]_{\lambda1} & i(xy) \ar[d]^{\lambda} \\ & xy } \]
commutes.
\end{proposition}

\proof
This amounts to commutativity of 
\[ \xymatrix @C3pc {
F(\myl_2) \ar[r]^-{\tilde{F}_{\myl_1,1,\myt_2}} \ar[d]_{F(\mylam_2)} & 
 F(\myl_1)\circ_1 F(\myt_2) \ar[d]^{F(\mylam_1)\circ_1 1} \\
 F(\myt_2) \ar@{=}[r]_{\tilde{F}_{\myt_1,1,\myt_2}} & F(\myt_1)\circ_1F(\myt_2) } \]
which follows by naturality of $\tilde{F}$.
\endproof

\begin{proposition}\label{prop:lr}
The composite 
\[ \xymatrix{ i \ar[r]^\rho & ii \ar[r]^{\lambda} & i } \]
is the identity.  
\end{proposition}

\proof
This amounts to commutativity of 
\[ \xymatrix{
F(\myl_0) \ar[r]^-{\tilde{F}_{\myl_1,1,\myl_0}}
\ar@{=}[dr]_{\tilde{F}_{\myt_1,1,\myl_0}} & F(\myl_1)\circ_1F(\myl_0)
\ar[d]^{F(\lambda_{1})\circ_{1}1} \\
& F(\myt_1)\circ_1F(\myl_0) } \]
in which the diagonal is an identity by counitality. The diagram
commutes by naturality of $\tilde{F}$ once again. 
\endproof

\begin{proposition}\label{prop:alr}
The composite 
\[ \xymatrix{ xy \ar[r]^{\rho1} & (xi)y \ar[r]^{\alpha} & x(iy) \ar[r]^{1\lambda} & xy } \]
is the identity.  
\end{proposition}

\proof
By naturality of $\tilde{F}$ once again, the diagram 
\[ \xymatrix{
F(\myt_2) \ar[d]_{\tilde{F}_{\myt_2,2,\myl_1}}
\ar@{=}[dr]^{\tilde{F}_{\myt_2,2,\myt_1}} \\
F(\myt_2) \circ_{2} F(\myl_1) \ar[r]_-{1 \circ_{2} F(\mylam_1)} &
F(\myt_2)\circ_{2}F(\myt_1) } \]
commutes, where the diagonal is an identity by counitality and the
horizontal is $1\lambda\colon x(iy)\to xy$. The vertical is $\alpha.\rho1$ by
Example~\ref{ex:a.r1}.
\endproof

\begin{proposition}\label{prop:pentagon}
  The pentagon
\[ \xymatrix{
((wx)y)z \ar[r]^{\alpha1} \ar[d]_{\alpha} & (w(xy))z \ar[r]^{\alpha} &
w((xy)z) \ar[d]^{1\alpha} \\
(wx)(yz) \ar[rr]_{\alpha} & & w(x(yz)) } \]
commutes.
\end{proposition}

\proof
By coassociativity, the diagram
\[ \xymatrix{
F(\myt_4) \ar[rr]^{\tilde{F}_{\myt_2,2,\myt_3}}
\ar[d]_{\tilde{F}_{\myt_3,3,\myt_2}} &&
F(\myt_2)\circ_2 F(\myt_3)
\ar[d]^{1 \circ_2\tilde{F}_{\myt_2,2,\myt_2}} \\
F(\myt_3)\circ_3 F(\myt_2) \ar[r]_-{\tilde{F}_{\myt_2,2,\myt_2}\circ_31 } & 
(F(\myt_2)\circ_2 F(\myt_2))\circ_3 F(\myt_2) \ar@{=}[r] & 
F(\myt_2)\circ_2(F(\myt_2)\circ_2 F(\myt_2)) } \]
commutes. The left vertical, right vertical, and lower path in this
diagram coincide with those in the statement of the proposition; and
the upper horizontal does too, by Example~\ref{ex:a.a1}.
\endproof

\subsection{The correspondence}

With these preparations, we are now ready to prove the following result.

\begin{theorem}
  If $w$ is any of the flavours \textnormal{$\ell$, c, ps, s} of weak
  morphism, the 2-functor $H^*\colon\Skeww\to\nc\cl\Algw$ of
  Corollary~\ref{cor:H}  is fully faithful,
  with image given by the \LNC-algebras. 
\end{theorem}

\proof
A skew monoidal category $A$, with structure corresponding to an
operad morphism $F\colon \cs\to\End(A)$, is sent to the normal colax
algebra $H^*(A)$ given by the composite $F\circ H\colon
\cl\to\End(A)$. This is indeed an \LNC-algebra by Example~\ref{ex:H}.

It is not hard to see that this is injective. The multiplication of
the skew monoidal category is encoded by $F(\ord 2,\ord
2,\bot_2)=FH(\myt_2)$, and the unit by $F(\ord
1,\emptyset,\bot_1)=FH(\myl_0)$, thus these are both determined by the
\LNC-algebra. Similarly the left unit map is given by $FH(\lambda_1)$,
the right unit map by $F\tilde{H}_{\myt_2,2,\myl_0}$, and the
associativity map by $F\tilde{H}_{\myt_2,2,\myt_2}$.

This proves that $H^*$ is injective on objects, and injectivity on
morphisms and 2-cells is similar but easier. 

Suppose conversely that $A$ is an \LNC-algebra, with corresponding
\LNC-morphism $F\colon \cl\to\End(A)$. 
Then we may define, as at the
beginning of the section, $m=F(\myt_2)\colon A^2\to A$, $i=F(\myl_0)\in
A$, and so on, and then by Propositions~\ref{prop:ar}, \ref{prop:al},
\ref{prop:lr}, \ref{prop:alr}, and \ref{prop:pentagon} this defines a
skew monoidal category. 
Furthermore, by Proposition~\ref{prop:spec-gen}, 
the resulting skew monoidal category is sent by $H^*$ to the original
\LNC-algebra. 

Now suppose that $f\colon A\to B$ is a lax morphism of
\LNC-algebras, and let $G\colon \cl\to\{f,f\}_\ell$ be the
corresponding colax morphism of operads. Since $d\circ G$ and
$c\circ G$ are \LNC-morphisms, and $d$ and $c$ are strict morphisms
which jointly reflect identities, it follows that $G$ is also an
\LNC-morphism. Now $G(\myt_2)$ has the form 
\[ \xymatrix{
A^2 \ar[r]^{f^2}_{~}="1" \ar[d]_{m} & B^2 \ar[d]^{m} \\
A \ar[r]_{f}^{~}="2" & B 
\ar@{=>}"1";"2" }
\]
while $G(\myl_0)$ has the form $i\to fi$. Furthermore, $\tilde{G}_{\myt_2,2,\myt_2}$ is
determined by $G(\myt_2)$ and the \LNC-algebra structures, and
encodes the associativity condition for $G(\myt_2)$. Similarly
$\tilde{G}_{\myl_2,2,\myl_0}$ and $G(\lambda_1)$ are  determined by the other data and
encode the unit conditions for $G(\myt_2)$. This proves fullness of
$H^*$ on lax morphisms. 

The cases of the other flavours of morphism and of 2-cells are
similar and left to the reader. 
\endproof

\subsection{The colax $L$-algebra associated to a skew monoidal category}\label{sec:colaxLalgebra}
We have  seen that skew monoidal categories correspond to
\LNC-algebras, and we have  given a concrete description of the skew monoidal category
associated to an \LNC-algebra. We now match this with a concrete description of
the \LNC-algebra associated to a skew monoidal category.  This will be used
in the companion paper  \cite{bourke-lack} to describe the \emph{skew multicategory} associated to a skew monoidal category.

Given a skew monoidal $C$, the corresponding $\cs$-algebra is specified
by an operad morphism $c\colon \cs \to \End(C)$, whose value at $x \in \cs_n$
we  write as $c(x)\colon C^n \to C$ and whose value at $\alpha\colon x
\to y$ we write as $c(\alpha)\colon c(x) \Rightarrow c(y)$.
At $(\ord 2, \ord 2,\bot) \in \cs_2$ and $(\ord 0,\ord 0,\bot) \in
\cs_0$ the corresponding functors are $\otimes\colon C^2 \to C$ and
$I\colon C^0 \to C$.  Every element of $\cs_n$ is obtained from the
above elements of $\cs_2$ and $\cs_0$ by operadic substitution; it
follows that the functors of the form $c(x)\colon C^n \to C$ are
precisely those obtained from $\otimes\colon C^2 \to C$ and $I\colon C^0 \to C$ by substitution.   

We note that the functors $c(x)\colon C^n \to C$ are 2-natural in
$C$\ -- more precisely, in the strict monoidal functors and monoidal
natural transformations of $\Skews$.  We mention this last abstract
point because we would like to say something about certain components
of the form $c(\alpha)\colon c(x) \to c(y)$ for $\alpha\colon x \to y
\in \cs_n$, whilst avoiding the syntax of $\cs_n$ itself.  To that
end, we point out that for each family of natural transformations
$\{\alpha_{C}\colon C \in \Skews\}$ natural in strict monoidal functors $F$ in the sense of the following diagram
\begin{equation}\label{eq:semantics}
\xy{
(0,0)*+{C^n}="00"; (40,0)*+{C}="10";
(0,-20)*+{D^n}="01";(40,-20)*+{\cd}="11";
{\ar@/^1pc/^{x(c)} "00"; "10"};{\ar@/_1pc/_{y(c)} "00"; "10"};
{\ar@/^1pc/^{x(d)} "01"; "11"};{\ar@/_1pc/_{y(d)} "01"; "11"};
{\ar_{F^{n}} "00"; "01"};{\ar^{F} "10"; "11"};
{\ar@{=>}^{\alpha_{C}}(20,3)*+{};(20,-3)*+{}};
{\ar@{=>}^{\alpha_{D}}(20,-17)*+{};(20,-23)*+{}};
}
\endxy
\end{equation}
there exists a \emph{unique} $\alpha\colon x \to y \in \cs_n$
with $c(\alpha) = \alpha_C$ for each skew monoidal $C$.  This is a
general syntax/semantics fact that holds for any (plain)  $\Cat$-operad.

Now the corresponding colax $\cl$-algebra $m\colon \cl \to \End(C)$ is given by the composite colax morphism of operads
 \begin{equation*}
 \xymatrix{
 \cl \ar[r]^-{H} & \cs \ar[r]^-{c} & \End(C)
 }
 \end{equation*}
 in which $H$ is the colax morphism of operads of Corollary~\ref{cor:H}.  This has components  
\begin{equation*}
 \xymatrix{
 \cl_n \ar[r]^-{H_{n}} & \cs_n \ar[r]^-{c_{n}} & [C^{n},C]
 }
 \end{equation*}
and substitution maps given by the natural transformations below.
\begin{equation*}
\xy{
(0,0)*+{\cl_n \times \cl_k}="00"; (40,0)*+{\cs_n \times \cs_k}="10";(80,0)*+{[C^{n},C] \times [C^{k},C]}="20";
(0,-20)*+{\cl_{n+k-1}}="01";(40,-20)*+{\cs_{n+k-1}}="11";(80,-20)*+{[C^{n+k-1},C]}="21";
{\ar_-{\circ_{i}} "00"; "01"};{\ar_-{H_{n+k-1}} "01"; "11"};{\ar^{H_{n} \times H_{k}} "00"; "10"};{\ar_{\circ_{i}} "10"; "11"};
{\ar@{=>}^{\tilde{H}}(15,-15)*+{};(25,-5)*+{}};{\ar^-{c_{n} \times c_{k}} "10"; "20"};{\ar_-{\circ_{i}} "20"; "21"};
{\ar_{c_{n+k-1}} "11"; "21"};
}
\endxy
\end{equation*}

Let us write $a_{1}\ldots a_{n}$ for the left bracketed $n$-fold
tensor product in $C$, so that $a_{1}\ldots a_{n}a_{n+1} =
(a_{1}\ldots a_{n}) \otimes a_{n+1}$.  By Theorem~\ref{thm:Q-left}
$H_{n}(l) = (\ord{n+1},\ord{n+1}/\{0\},\bot_{n+1})$ is the
\emph{initial object} of $\cs_n$, consisting of the ordinal
$\ord{n+1}$ with $0$ omitted, and $\bot_{n+1} \in T_{n+1}$ the least
element of the Tamari lattice, corresponding to the left bracketing of
$n+1$-elements.  Accordingly $m_{l_{n}}=c_{n} \circ H_{n}(l)\colon C^{n} \to C$ has value $$m_{l_{n}}(a_{1},\ldots ,a_{n})  = ia_{1}\ldots a_{n}$$
the left bracketing of the $n+1$-tuple $(i,a_{1},\ldots,a_{n})$.

By Theorem~\ref{thm:Q-left} we also have $H_{n}(t) =
(\ord{n},\ord{n},\bot_{n})$.  It follows that $$m_{t_{n}}(a_{1},\ldots
,a_{n})=a_{1}\ldots a_{n}$$ the leftmost bracketing.  At
$\lambda\colon l \to t \in \cl_n$ for $n>0$ the induced map from
$m_{l_{n}}(\overline{a}) \to m_{t_{n}}(\overline{a})$ has component
$\lambda a_{1}\ldots a_{n}\colon ia_{1}\ldots a_{n} \to a_{1}\ldots a_{n}$.

With regards substitution, it follows from
Remark~\ref{remark:injection} that for any $(x,y) \in \cl_n \times \cl_k$
the morphism $\tilde{H}_{x,y}\colon H_{n+k-1}(x \circ_{i} y) \to
H_{n}(x) \circ_{i} H_{k}(y) \in \cs_{n+k-1}$ is \emph{unique} and,
moreover, an \emph{$\Fsk$-injection}.  Uniqueness allows us to say
that the substitution component $m_{x \circ_{i} y} \to m_{x} \circ_{i}
m_{y}$ is the unique natural transformation that exists for all skew
monoidal categories $C$ and is natural in the sense of
Diagram~\eqref{eq:semantics}.  
Furthermore, by Proposition~\ref{prop:Fsk-inj-fact1} each
$\Fsk$-injection admits a canonical factorisation as a swell morphism
(corresponding to applications of $\rho$) followed by a Tamari
morphism (corresponding to applications of $\alpha$).  Accordingly
each substitution is obtained by repeated applications of the right
unit maps $\rho$ followed by applications of associativity maps $\alpha$, each possibly tensored on either side.  
For instance $m_{l}(a,b,c,d) \to m_{l}(m_{t}(a,b),m_{l}(c,d))$ is the map given by
\begin{equation*}
\xymatrix @C1.4pc {
(((ia)b)c)d \ar[rr]^-{((i(ab))\rho)d} && (((ia)b)(ic))d  \ar[rr]^{((\alpha(ic))d } && (((i(ab))(ic))d \ar[r]^-{\alpha} & (i(ab))((ic)d).
}
\end{equation*}


\bibliographystyle{plain}

\end{document}